%% file: SurfacesTexArxivRevised.tex
\documentclass[11pt]{article}
\usepackage{placeins}
\usepackage{listings}
\usepackage[utf8]{inputenc} 
\usepackage{amsmath}
\usepackage{amsfonts}
\usepackage{calligra}
\usepackage{geometry} 
\usepackage{graphicx} 
\usepackage{mathpazo}
\usepackage{tikz}
\usepackage{emptypage}
\usepackage{cite}
\usepackage{subcaption}
\usepackage{authblk}
\usepackage{lineno}
\usepackage[
    type={CC},
    modifier={by-nc-nd},
    version={4.0},
]{doclicense}
\usetikzlibrary{shapes.geometric}
\usetikzlibrary{arrows.meta}
\usetikzlibrary{perspective}
\graphicspath{{figures/}}
\usepackage{booktabs} 
\usepackage{array} 
\usepackage{paralist} 
\usepackage{verbatim} 
\usepackage{bm}
\usepackage{float}
\usepackage{fancyhdr} 
\usepackage{tcolorbox}
\usepackage{csvsimple}

\usepackage{import}
\usepackage{xifthen}
\usepackage{pdfpages}
\usepackage{transparent}

\numberwithin{equation}{section}

\usepackage[nottoc,notlof,notlot]{tocbibind} 
\usepackage[titles,subfigure]{tocloft} 
\usepackage{tcolorbox}
\newcommand{\R}{\mathbb{R}}

\newcommand{\pp}[2]{\frac{\partial {#1}}{\partial {#2}}}

\renewcommand{\vec}[1]{\mathbf{#1}}
\newcommand{\vecg}[1]{\boldsymbol{#1}}
\newcommand{\xf}{\vec x_f}
\newcommand{\magn}[1]{\lvert {#1} \rvert}

\makeatletter
\newcommand{\Rmnum}[1]{\expandafter\@slowromancap\romannumeral #1@}
\newcommand{\whvec}[1]{\widehat{\vec{#1}}}
\newcommand{\whvecg}[1]{\widehat{\vecg{#1}}}
\newcommand{\bell}{\bm{\hat{\ell}}}

\usepackage{algorithm}
\usepackage{algpseudocode}
\usepackage{enumitem}

\usepackage{scalerel,stackengine}
\usepackage{color}   
\usepackage{hyperref}
\hypersetup{
    colorlinks=true,
    linktoc=all,     
    linkcolor=blue,  
}
\stackMath
\newcommand\widecheck[1]{%
\savestack{\tmpbox}{\stretchto{%
  \scaleto{%
    \scalerel*[\widthof{\ensuremath{#1}}]{\kern-.6pt\bigwedge\kern-.6pt}%
    {\rule[-\textheight/2]{1ex}{\textheight}}
  }{\textheight}%
}{0.5ex}}%
\stackon[1pt]{#1}{\scalebox{-1}{\tmpbox}}%
}

\newcolumntype{?}{!{\vrule width 1pt}}

\begin{document}
\author[1,\footnote{\noindent Corresponding author: \textit{dferranti@wpi.edu},  Permanent address: Department of Mathematical Sciences, Worcester Polytechnic Institute, Worcester, MA}]{Dana Ferranti}
\author[1]{Ricardo Cortez}
\affil[1]{Department of Mathematics, Tulane University, New Orleans, LA}
\makeatother

\title{Regularized Stokeslet Surfaces}
\date{}
\maketitle

\vspace{-2em}
\begin{abstract}
A variation of the Method of Regularized Stokeslets (MRS) in three dimensions is developed for triangulated surfaces with a piecewise linear force density. The work extends the regularized Stokeslet segment methodology used for piecewise linear curves. By using analytic integration of the regularized Stokeslet kernel over the triangles, the regularization parameter $\epsilon$ is effectively decoupled from the spatial discretization of the surface. This is in contrast to the usual implementation of the method in which the regularization parameter is chosen for accuracy reasons to be about the same size as the spatial discretization. The validity of the method is demonstrated through several examples, including the flow around a rigidly translating/rotating sphere and in the squirmer model for ciliate self-propulsion. Notably, second order convergence in the spatial discretization for fixed $\epsilon$ is demonstrated. Considerations of mesh design and choice of regularization parameter are discussed, and the performance of the method is compared with existing quadrature-based implementations.
\end{abstract}


\section{Introduction}
The dynamics of viscous-dominated fluids are relevant in a variety of scientific investigations. Often, the subjects of interest are at the microscale and intersect with biology, like bacteria or sperm motility \cite{grayhancock,fauciDillon,gaffney11,schuech}, and ciliary propulsion \cite{brennen77,ding2014,chakra22}. Other topics include suspension dynamics \cite{suspensionmorris}, phoretic particles \cite{phoretic1,phoretic2}, and biomedical technology \cite{microrobots,magdanz2020ironsperm,zhou2021}.

In the viscous-dominated regime of Newtonian fluids, fluid motion is described by the Stokes equations,

\begin{equation}
  \label{eq:stokes}
  \begin{split}
    \vec 0 &= -\nabla p + \mu \Delta \vec u + \vec F\\
    \nabla \cdot \vec u &=0
  \end{split}
\end{equation}where $p$ is the dynamic pressure of the fluid, $\mu$ is its dynamic viscosity, $\vec u$ is the velocity, and $\vec F$ is an external force per unit volume. For simulating flows generated by external forces, a popular approach introduced by Cortez and extended with collaborators is the method of regularized Stokeslets (MRS) \cite{cortez2001,cortez2005}. The method works by spatially spreading a force at the point of application through a blob function, $\phi_{\epsilon}$, where $\epsilon>0$ is a parameter controlling the size of the spreading. The blob function $\phi_{\epsilon}\left(\vec x \right)$ is a smooth function that satisfies $\int \int \int_{\R^3} \phi_{\epsilon}\left (\vec x \right) \ dV(\vec x)=1$ and can be thought of as an approximation of the Dirac delta distribution.

For a single regularized force applied at an arbitrary point $\vec y$ in a fluid of infinite expanse in three dimensions, the Stokes equations \eqref{eq:stokes} become

\begin{equation}
  \label{eq:regStokes}
  \begin{split}
    \vec 0 &= -\nabla p(\vec x) + \mu \Delta \vec u(\vec x) + \vec f \phi_{\epsilon}(\vec x-\vec y)\\
    \nabla \cdot \vec u(\vec x) &=0
  \end{split}
\end{equation}The regularized velocity field that solves \eqref{eq:regStokes} can be written in the form

\begin{equation}
  \label{eq:regStokesSoln}
  \left (8 \pi \mu \right )u(\vec x) = \vec S^{\epsilon}(\vec x-\vec y)\cdot \vec f
\end{equation}where $\vec S^{\epsilon}$ is called the regularized Stokeslet, or regularized Stokeslet kernel, whose particular form depends on the blob function $\phi_{\epsilon}$. We include the factor ${8\pi \mu}$ to agree with the common form of the analogous singular Stokeslet \cite{pozr1992}. The MRS is often used for modeling bodies (e.g. organisms, particles) which exert forces on the surrounding fluid. In a three-dimensional fluid, the MRS framework allows for the representation of these bodies as one-dimensional curves (line integral of regularizd Stokeslets) or two-dimensional surfaces (surface integral of regularized Stokeslets.) Scientific work that has used the MRS for this purpose include investigations of cilia and flagellar motion in Stokes and Brinkman flows \cite{olson13,olson16,brinkman19}, flows due to phoretic particles, \cite{montenegro15,varma2018}, and detailed models of microorganisms \cite{nguyen2019,lim2019}. The method has also been extended for use for use in bounded geometries, including the half-plane above a no-slip wall \cite{ainley,cortezvarela}, outside a solid sphere \cite{wrobel2016}, and periodic domains in two and three dimensional flow \cite{mannan,hoffmann2,hoffmann1}.

In practice, when using the MRS in the context of surface integrals solved via quadrature, a choice must be made in balancing the spatial discretization of the surface, $h$, with the regularization parameter, $\epsilon$. The complicating factor is the ``near-singularity'' of the regularized Stokeslet kernel for small $\epsilon$ compared to $h$. Ideally, one would like to use a small $\epsilon$ so as to minimize the error introduced by the regularization. However, the smaller $\epsilon$ is chosen, the more ``nearly singular'' the regularized Stokeslet kernel becomes. Even if a high order quadrature is used, there is no escaping the fact that the error estimates will depend inversely on a power of $\epsilon$ \cite{cortez2005}.

Several methods have been introduced to ameliorate this dependence of the spatial discretization on the regularization parameter. A boundary-element regularized Stokeslet method was proposed by Smith \cite{smith2009} in which the force discretization is decoupled from the kernel evaluations. In the simplest implementation of the method, the forces are constant over each mesh element and a high order quadrature method is used to evaluate the regularized Stokeslet kernel over the mesh. The key insight is that the set of points used in the discretization of the force density can be decoupled from the points used for the regularized Stokeslet kernel evaluations. A high order quadrature method can then be used without significantly increasing the computational cost of the method.

Another method proposed by Barrero-Gil addressed the problem by using ``auxiliary Stokeslets'' \cite{bgil}. In essence, the method evaluates the contribution to the velocity from forces in the near-field of an evaluation point by introducing more Stokeslets (the ``auxiliary'' ones) in a small region which are then spatially averaged. This extra computational effort is done when the force point and the evaluation point are the same, but not when they are different. With such selective refinement, better accuracy is achieved without significant additional cost.

We note one last method, introduced by Smith, called the nearest-neighbor interpolation regularized Stokeslet method \cite{smith18,gallagher}. This method uses two sets of points (potentially overlapping): a coarse set of force points and a finer set of quadrature points. The quadrature points are used for the quadrature evaluation and take on the forces of their nearest neighbors in the set of force points. The quadrature discretization can be made much finer than the force discretization without significantly increasing the cost of the method, since while the number of kernel evaluations may be large if there are many quadrature points, the number of degrees of freedom is proportional to the number of force points.

All of the above methods have made progress in weakening the $h,\epsilon$ dependence. This work takes a different approach. The method presented herein requires using a triangulated mesh of a surface. On each triangle, a linear force distribution is generated based on the values at the vertices, which results in a continuous piecewise linear force distribution over the entire surface. The fluid velocity due to this piecewise linear force density is evaluated analytically based on expressions we derive in the next section. Since there is no quadrature, the error is due to the regularization, the imposition of the piecewise linear force density, and the approximation of the surface by flat triangles. In effect, the regularization can be made much smaller than the size of the force discretization without any tuning. In the following section, we derive the method in detail. The second part of the paper discusses validation examples and applications.

\section{Derivation of the Method}
\begin{figure}[ht!]
  \centering
  \begin{subfigure}{0.49\textwidth}
    \centering
\begin{tikzpicture}
  \draw[blue,arrows={-Stealth}] (0,0) -- (-1,1) node[anchor=east]{$\vec f_0$};
  \draw[blue,arrows={-Stealth}] (3,0) -- (3.5,0.75) node[anchor=west]{$\vec f_1$};
  \draw[blue,arrows={-Stealth}] (1,2) -- (1.1,3) node[anchor=west]{$\vec f_2$};
  \fill[draw=black]
     (0,0) circle (3pt) node[anchor=north] {$\vec y_0$}
  -- (3,0) circle (3pt) node[anchor=north] {$\vec y_1$}
  -- (1,2) circle (3pt) node[anchor=east] {$\vec y_2$}
  -- (0,0);
  \draw[black, dashed, arrows = {-Latex[width=6pt,length=6pt]}] (0,0) -- (-1,2.4)  node[anchor=south]{};
  \fill[draw=black](-0.5,1.25) node[anchor=west] {$\vec x_0$};
  \fill[draw=black](-1,2.5) circle (3pt) node[anchor=south] {$\vec x_f$};
  \draw[magenta ,arrows={-Stealth}] (0.1,0) -- (2,0) node[anchor=north east]{\scriptsize $-L_1\alpha \whvec v $};
  \draw[blue ,loosely dashed, thick] (0.1,0) -- (2,0);
  \draw[magenta ,arrows={-Stealth}] (2,0) -- (1.15,0.85); 
  \draw[blue ,loosely dashed, thick] (2,0) -- (1.15,0.85);
   \draw (0.9,0.7) node[below,magenta] {\scriptsize $-L_2 \beta \whvec w$};
   \fill[draw=magenta,dashed] (1.25,0.75) -- (2,2) node[anchor=west,color=magenta]{\scriptsize $\vec y(\alpha,\beta)=\vec y_0-L_1\alpha \whvec v - L_2\beta \whvec w$};   
   \fill[draw=blue,dashed] (1.25,0.75) -- (1.5,2.4) node[anchor=west,color=blue]{\scriptsize $\vec f(\alpha,\beta)=\vec f_0 + \alpha (\vec f_1 - \vec f_0) + \beta (\vec f_2 - \vec f_1)$};
 \end{tikzpicture}
 \caption{Triangle region $\Omega$ in physical space. \label{fig:omega}}
\end{subfigure} 
\begin{subfigure}{0.5\textwidth}
  \centering
\begin{tikzpicture}
  \fill[draw=black]
     (0,0) circle (2pt) node[anchor=north] {\scriptsize $(0,0)$}
  -- (2,0) circle (2pt) node[anchor=north] {\scriptsize $(1,0)$}
  -- (2,2) circle (2pt) node[anchor=east] {\scriptsize $(1,1)$}
  -- (0,0);
  \draw[magenta ,arrows={-Stealth}] (0.1,0) -- (1.5,0) node[anchor=north east]{\scriptsize $-\alpha \whvec e_1 $};
  \draw[magenta ,arrows={-Stealth}] (1.5,0) -- (1.5,0.7);
  \draw (1.5,0.8) node[below left, magenta]{\scriptsize $-\beta \whvec e_2$};
  \draw[magenta, arrows={-Stealth}] (1.9,1.9) -- (1.3,1.3) node[anchor = west]{\scriptsize $-\whvec d$};
  \fill[draw=black](1,-0.1) node[anchor=south]{\scriptsize $S_1$};
  \fill[draw=black](2,1) node[anchor=west]{\scriptsize $S_2$};
  \fill[draw=black](1,1) node[anchor=south east]{\scriptsize $S_3$};
\end{tikzpicture}
\caption{Triangle region $\Omega_p$ in parameter space. \label{fig:omegap}}
\end{subfigure}
\caption{A sketch of the triangle in physical space and in the parameter space.}
\label{fig:triangle}
\end{figure}
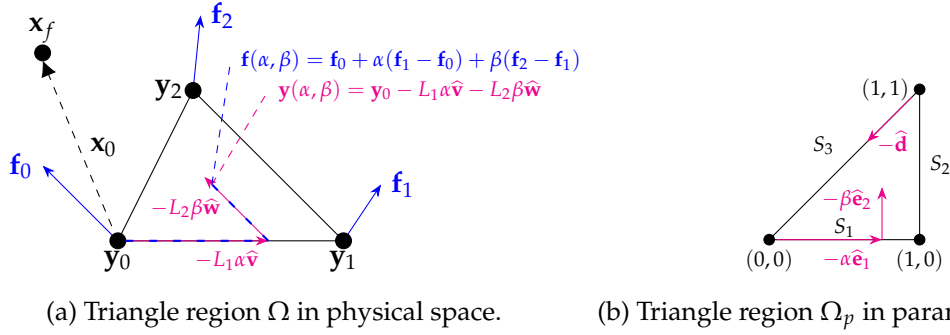

We begin by introducing the essential notation used throughout. The basic building block for the method is a triangle. A sketch accompanying the following description is shown in Figure \ref{fig:omega}.

The three vertices of the triangle are denoted $\vec y_0, \ \vec y_1, \ \vec y_2$. The unit vector $\whvec{v}$ lies in the direction of $\vec y_0 - \vec y_1$ and  $\whvec w$ is the unit vector that lies in the direction of $\vec y_1 - \vec y_2$. The distances $\lvert \vec y_0 - \vec y_1 \rvert$ and $\lvert \vec y_1 - \vec y_2 \rvert$ are $L_1$ and $L_2$, respectively.

A point in or on the boundary of the trianglular region $\Omega$ is written parametrically as $\vec y(\alpha, \beta) = \vec y_0 - \alpha L_1 \whvec v - \beta L_2 \whvec w$, where $\alpha, \beta$ are dimensionless parameters such that $0 \leq \alpha \leq 1$, and $\beta = \beta(\alpha)$ such that for fixed $\alpha$, $0 \leq \beta(\alpha) \leq \alpha$. The force density $\vec f(\alpha,\beta)$ is a linear interpolation of the force densities $\vec f_0, \ \vec f_1, \ \vec f_2$, at the respective points $\vec y_0, \ \vec y_1, \ \vec y_2$, i.e. $\vec f(\alpha,\beta)=\vec f_0 + \alpha (\vec f_1-\vec f_0) + \beta(\vec f_2-\vec f_1)$. For conciseness, we will write this as $\vec f(\alpha,\beta)=\vec f_0 + \alpha \vec f_a + \beta \vec f_b$, where $\vec f_a = \vec f_1 - \vec f_0$ and $\vec f_b = \vec f_2 - \vec f_1$.

The triangle parameterization in $(\alpha,\beta)$ implicitly maps the triangular region in physical space, $\Omega$, to the triangular region in parameter space, $\Omega_p$, shown in Figure \ref{fig:omegap}. This mapping is linear and its Jacobian is twice the area of the triangle, $BH$ (read ``base $\times$ height''.)

The net force $\vec F$ and torque $\vec M$ (with respect to a central point $\vec y_c$) over the interior of a single triangle, $\partial \Omega$, are respectively,

\[ \vec F = \int \int_{\Omega} \vec f(\vec y) \ d \vec y = (BH) \int \int_{\Omega_p} \vec f(\vec y(\alpha,\beta)) d\beta \ d\alpha = (BH) \int_0^1 \int_0^{\alpha} \vec f_0 + \alpha \vec f_a + \beta \vec f_b \ d \beta \ d \alpha \]and

\begin{equation*}
  \begin{split}
   \vec M &=  \int \int_{\Omega} (\vec y - \vec y_c) \times \vec f(\vec y) \ d \vec y = (BH) \int \int_{\Omega_p} (\vec y(\alpha,\beta) - \vec y_c) \times \vec f(\vec y(\alpha,\beta)) \ d\beta \ d\alpha \\
          &= (BH)\int_0^1 \int_0^{\alpha}  \left(\vec y_0 + \alpha (\vec y_1-\vec y_0) + \beta (\vec y_2-\vec y_1) - \vec y_c \right ) \times \left(\vec f_0 + \alpha \vec f_a + \beta \vec f_b \right) \ d \beta \ d \alpha 
  \end{split}
\end{equation*}After substituting $\vec f_a = \vec f_1 - \vec f_0$ and $\vec f_b = \vec f_2 - \vec f_1$ and simplifying, these expressions become 

\begin{align}
  \vec F &= \frac{BH}{6}\Big (\vec{f_0} + \vec f_1 + \vec f_2 \Big ) \label{eq:netForce}\\
  \vec M &= \frac{BH}{24} \Big ((2\vec{y_0} + \vec y_1 + \vec y_2 - 4 \vec y_c) \times \vec f_0 + (\vec y_0 + 2 \vec y_1 + \vec y_2 - 4 \vec y_c) \times \vec f_1 \nonumber \\
   &{} \hspace{4em} (\vec y_0 + \vec y_1 + 2 \vec y_2 - 4 \vec y_c) \times \vec f_2 \Big)   \label{eq:netTorque}
\end{align}

We are interested in analytically solving for the velocity and pressure at a field point $\xf$ due to the force distribution $\vec f(\vec y)$ (force per unit area) over the triangle. More precisely, we wish to solve

\begin{equation}
  \label{eq:stokessurfeq}
  \begin{split}
    \vec 0 &= -\nabla p(\xf) + \mu \Delta \vec u(\xf) + \int \int_{\Omega} \vec f(\vec y) \phi_{\epsilon}(\magn{\xf-\vec y}) \ d \vec y \\
    \nabla \cdot \vec u(\xf) &=0 
  \end{split}
\end{equation}where the blob function $\phi_{\epsilon}(r)$ is

\begin{equation}
  \label{eq:blob}
  \phi_{\epsilon}(r)= \frac{15 \epsilon^4}{8\pi(r^2+\epsilon^2)^{7/2}}
\end{equation}By the superposition principle, the solution to \eqref{eq:stokessurfeq} is a surface integral of regularized Stokeslets multiplying the force density, 
\begin{equation}
  \label{eq:stokessurfsol}
  \vec u (\xf) = \frac{1}{8\pi\mu}\int \int_{\Omega} \vec S^{\epsilon}(\xf, \vec y) \cdot \vec f(\vec y) \ d \vec y
\end{equation}where $\vec S^{\epsilon}(\vec x, \vec y)$ is the particular regularized Stokeslet kernel for the blob function in \eqref{eq:blob}, given in component form as, 

\begin{equation}
  \label{eq:regstokes}
  \begin{split}
    S_{ij}^{\epsilon}(\vec x, \vec y) &= \left( \frac{1}{R}+\frac{\epsilon^2}{R^3}\right)\delta_{ij} + \frac{(x_i - y_i)(x_j - y_j)}{R^3} \\
   \text{where } R^2(\vec x, \vec y):&= \lvert \vec x - \vec y \rvert^2 + \epsilon^2
  \end{split}
\end{equation}When the context is clear, we will omit the arguments of $R(\vec x, \vec y)$ and write it simply as $R$. In terms of the ($\alpha,\beta$) parameterization, the velocity $\vec u(\vec x_f)$ is written,

\begin{equation}
  \label{eq:stokessurfsolparam}
  \vec u (\xf) = \frac{BH}{8 \pi \mu} \int \int_{\Omega_p} \vec S^{\epsilon}(\xf, \vec y(\alpha,\beta)) \cdot \vec f(\alpha,\beta)) \ d \beta \ d\alpha
\end{equation}

Substituting $\vec y(\alpha,\beta)$ for $\vec y$ in \eqref{eq:regstokes}, we note that this parameterized form of the regularized Stokeslet kernel $\vec S^{\epsilon}$ has terms quadratic in $\alpha,\beta$. Since the force density $\vec f(\alpha,\beta)$ is linear in $\alpha,\beta$, the integrand $\vec S^{\epsilon}(\xf, \vec y(\alpha,\beta))\cdot \vec f(\alpha,\beta)$ is a cubic polynomial in $\alpha,\beta$ with constant vector coefficients,

\begin{equation}
  \label{eq:expansion}
  \begin{split}
    \vec S^{\epsilon}(\xf,\vec y(\alpha,\beta))\cdot \vec f(\alpha,\beta) &= \frac{\vec f_0}{R}  + \frac{\vec P_{0,0}}{R^3} + \alpha \left (\frac{\vec f_a}{R}+\frac{\vec P_{1,0}}{R^3} \right ) + \beta \left (\frac{\vec f_b}{R}+ \frac{\vec P_{0,1}}{R^3}\right ) \\
    &+ \alpha^2\frac{\vec P_{2,0}}{R^3} + \alpha \beta \frac{\vec P_{1,1}}{R^3} + \beta^2 \frac{\vec P_{0,2}}{R^3} \\
    &+ \alpha^3 \frac{\vec P_{3,0}}{R^3} + \alpha^2 \beta \frac{\vec P_{2,1}}{R^3}+ \alpha \beta^2 \frac{\vec P_{1,2}}{R^3} + \beta^3 \frac{\vec P_{0,3}}{R^3}
  \end{split}
\end{equation}where the $\vec P_{i,j}$ are the constant vector coefficients listed in Table \ref{Tab:Coefficients}.

Note that \eqref{eq:expansion} includes only terms of the form $\alpha^m \beta^n R^{-q}$ where $m+n \leq q$ and $q=\{1,3\}$. Exact integration of \eqref{eq:stokessurfsolparam} then requires knowing how to compute (for the aforementioned values of $m,n,q$),

\begin{equation}
  \label{eq:Tmnq}
  T_{m,n,q}:= \int \int_{\Omega_p}\alpha^m \beta^n R^{-q} \ d \ \vec s(\alpha,\beta)
\end{equation}We write the final formula for the velocity, $\vec u(\xf)$, concisely as 

\begin{equation}
  \label{eq:velotri}
  \begin{split}
    \frac{8 \pi \mu}{BH}\vec u(\xf) &= \vec f_0 \ T_{0,0,1} + \vec P_{0,0} \ T_{0,0,3} + \vec f_a \ T_{1,0,1} + \vec P_{1,0} \ T_{1,0,3} + \vec f_b \ T_{0,1,1}\\
                                    &+ \vec P_{0,1} \ T_{0,1,3} + \vec P_{2,0} \ T_{2,0,3} + \vec P_{1,1} \ T_{1,1,3} + \vec P_{0,2} \ T_{0,2,3} \\
                                    &+ \vec P_{3,0} \  T_{3,0,3} + \vec P_{2,1} \ T_{2,1,3} + \vec P_{1,2} \ T_{1,2,3} + \vec P_{0,3} \ T_{0,3,3}
  \end{split}
\end{equation}

\begin{table}[ht]
\centering
\begin{tabular}{c|l}%
  \hline
  $  \vec P_{0,0}$ & $\epsilon^2 \vec f_0 + (\vec f_0 \cdot \vec x_0) \vec x_0$ \\
  \hline
  $\vec P_{1,0}$ & $\epsilon^2 \vec f_a + (L_1 \vec f_0 \cdot \whvec v + \vec f_a \cdot \vec x_0 )\vec x_0 + (L_1 \vec f_0 \cdot \vec x_0) \whvec v$   \\
  \hline
  $\vec P_{0,1}$ & $\epsilon^2 \vec f_b + (L_2 \vec f_0 \cdot \whvec w + \vec f_b \cdot \vec x_0 )\vec x_0 + (L_2 \vec f_0 \cdot \vec x_0) \whvec w$   \\
  \hline
  $\vec P_{2,0}$ & $(L_1 \vec f_a \cdot \whvec v) \vec x_0 + (L_1^2 \vec f_0 \cdot \whvec v + L_1 \vec f_a \cdot \vec x_0 ) \whvec v$  \\
  \hline
  $\vec P_{1,1}$ & \begin{tabular}{@{}l@{}}$(L_1 \vec f_b \cdot \whvec v + L_2 \vec f_a \cdot \whvec w)\vec x_0 + (L_1L_2 \vec f_0 \cdot \whvec w + L_1 \vec f_b \cdot \vec x_0 )\whvec v$ \\ $+ (L_1 L_2 \vec f_0 \cdot \whvec v + L_2 \vec f_a \cdot \vec x_0) \whvec w$ \end{tabular}\\
  \hline
  $\vec P_{0,2}$& $(L_2 \vec f_b \cdot \whvec w) \vec x_0 + (L_2^2 \vec f_0 \cdot \whvec w + L_2 \vec f_b \cdot \vec x_0 ) \whvec w$ \\
  \hline
  $\vec P_{3,0}$& $(L_1^2 \vec f_a \cdot \whvec v)\whvec v$ \\
  \hline
  $\vec P_{2,1}$& $(L_1 L_2 \vec f_a \cdot \whvec w + L_1^2 \vec f_b \cdot \whvec v) \whvec v + (L_1 L_2 \vec f_a \cdot \whvec v) \whvec w$ \\
  \hline
  $\vec P_{1,2}$& $(L_1 L_2 \vec f_b \cdot \whvec v + L_2^2 \vec f_a \cdot \whvec w) \whvec w + (L_1 L_2 \vec f_b \cdot \whvec w) \whvec v$ \\
  \hline
  $\vec P_{0,3}$& $(L_2^2 \vec f_b \cdot \whvec w)\whvec w$ \\ \hline
\end{tabular}
\caption[Vector coefficients in equation \eqref{eq:expansion}]{Vector coefficients in equation \eqref{eq:expansion}. As in Figure \ref{fig:triangle}, we write $\vec x_0=\vec x_f-\vec y_0$.}
\label{Tab:Coefficients}
\end{table}

\subsection{Recursion Formulas}
A recurrence formula relating the $T_{m,n,q}$ can be established which simplifies the work in evaluating \eqref{eq:velotri} for all of the cases. This is analagous to the recurrence formula established for regularized Stokeslet segments \cite{cortez18}. In the following, we refer to the triangle shown in Figure \ref{fig:omegap} whose boundary is $\partial \Omega_p$. This is a right triangle with orthogonal sides in the directions $\whvec e_1$, $\whvec e_2$, and diagonal $\whvec d = \frac{1}{\sqrt{2}}\left ( \whvec e_1 + \whvec e_2 \right)$.

\begin{figure}[ht]
    \centering
    \scalebox{0.7}{
    \def\svgwidth{\columnwidth}
    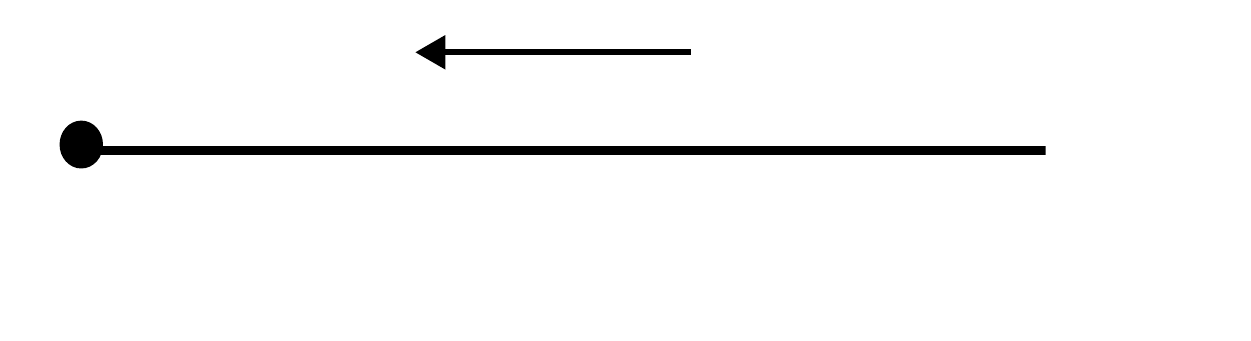
}
    \caption[Geometry of line segment integrals.]{Geometry of line segment integrals \ref{eq:Slm1} and \ref{eq:Slp1}. The length of the line segment is $L$ and the $\bell$ points from $\vec y(1)$ to $\vec y(0)$.}
    \label{fig:segment}
  \end{figure}

We begin by using integration by parts in $\alpha$ and $\beta$ to establish the respective identities,

\begin{align}
  &\int \int_{\Omega_p} \pp{}{\alpha} \left (\alpha^m \beta^n R^{-(q-2)} \right ) \ d \vec s(\alpha,\beta) = \int_{\partial \Omega_p} \alpha^m(\theta) \beta^n(\theta) \ R^{-(q-2)} \ \whvec{n}(\theta) \cdot \whvec e_1 \ d\theta =: A_{m,n,q-2} \label{eq:intbypartsa}\\
  &\int \int_{\Omega_p} \pp{}{\beta} \left (\alpha^m \beta^n R^{-(q-2)} \right ) \ d \vec s(\alpha,\beta) = \int_{\partial \Omega_p} \alpha^m(\theta) \beta^n(\theta) \ R^{-(q-2)} \ \whvec{n}(\theta) \cdot \whvec e_2 \ d\theta =: B_{m,n,q-2}\label{eq:intbypartsb}
\end{align}where $\theta$ parameterizes the triangle boundary and $\whvec n(\theta)$ is the outward pointing normal to the side of the triangle at the point corresponding to $\theta$. One can also directly differentiate the integrand of the integral on the left hand side to get a different pair of identities,

\begin{equation}
  \begin{split}
  &\int \int_{\Omega_p} \pp{}{\alpha} \left (\alpha^m \beta^n R^{-(q-2)} \right ) \ d \vec s(\alpha,\beta)\\
  &= m T_{m-1,n,q-2}-(q-2) \left(L_1(\vec x_0 \cdot \whvec v) T_{m,n,q} + L_1^2T_{m+1,n,q}+L_1L_2 (\whvec v \cdot \whvec w) T_{m,n+1,q} \right)
  \label{eq:diffida}
\end{split}
\end{equation}
\begin{equation}
  \begin{split}
  &\int \int_{\Omega_p} \pp{}{\beta} \left (\alpha^m \beta^n R^{-(q-2)} \right ) \ d \vec s(\alpha,\beta)\\
  &= n T_{m,n-1,q-2}-(q-2) \left(L_2(\vec x_0 \cdot \whvec w) T_{m,n,q} + L_2^2T_{m,n+1,q}+L_1L_2 (\whvec v \cdot \whvec w) T_{m+1,n,q} \right)
  \label{eq:diffidb}
\end{split}
\end{equation}Equating \eqref{eq:intbypartsa} with \eqref{eq:diffida} and \eqref{eq:intbypartsb} with \eqref{eq:diffidb}, we have

  \begin{align}
    A_{m,n,q-2}&= m T_{m-1,n,q-2} \nonumber  \\
    &\hspace{1em} -(q-2) \left(L_1(\vec x_0 \cdot \whvec v) T_{m,n,q} + L_1^2T_{m+1,n,q}+L_1L_2 (\whvec v \cdot \whvec w) T_{m,n+1,q} \right) \label{eq:ida}\\
    B_{m,n,q-2}&= n T_{m,n-1,q-2} \nonumber \\
    &\hspace{1em} -(q-2) \left(L_2(\vec x_0 \cdot \whvec w) T_{m,n,q} + L_2^2T_{m,n+1,q}+L_1L_2 (\whvec v \cdot \whvec w) T_{m+1,n,q} \right) \label{eq:idb}
  \end{align}A recursion for the first index is established by subtracting $\frac{L_1}{L_2}(\whvec v \cdot \whvec w) \times$ \eqref{eq:idb} from \eqref{eq:ida} and rearranging to arrive at,

  \begin{equation}
    \begin{split}
      T_{m+1,n,q} = \frac{1}{(\whvec v \cdot \whvec w)^2-1} \Bigg ( &\frac{A_{m,n,q-2}}{L_1^2(q-2)} - (\whvec v \cdot \whvec w)\frac{B_{m,n,q-2}}{L_1L_2(q-2)} -\frac{m}{L_1^2(q-2)}T_{m-1,n,q-2} \\
                                                                    &+(\whvec v \cdot \whvec w)\frac{n}{L_1L_2(q-2)}T_{m,n-1,q-2} \\
                                                                    & +\frac{\vec x_0 \cdot \whvec v - (\whvec v \cdot \whvec w) \vec x_0 \cdot \whvec w}{L_1}T_{m,n,q} \Bigg ), \ q\neq 2
    \end{split}
    \label{eq:Tmp1}
  \end{equation}A recursion can be similarly found for the second index by subtracting $\frac{L_2}{L_1} \times$ \eqref{eq:ida} from \eqref{eq:idb} and rearranging,
  \begin{equation}
    \begin{split}
      T_{m,n+1,q} = \frac{1}{(\whvec v \cdot \whvec w)^2-1} \Bigg ( &\frac{B_{m,n,q-2}}{L_2^2(q-2)} - (\whvec v \cdot \whvec w)\frac{A_{m,n,q-2}}{L_1L_2(q-2)}-\frac{n}{L_2^2(q-2)}T_{m,n-1,q-2}\\
                                                                    &+(\whvec v \cdot \whvec w)\frac{m}{L_1L_2(q-2)}T_{m-1,n,q-2}\\
                                                                    &+ \frac{\vec x_0 \cdot \whvec w - (\whvec v \cdot \whvec w) \vec x_0 \cdot \whvec v}{L_2}T_{m,n,q} \Bigg ), \ q\neq 2
    \end{split}
    \label{eq:Tnp1}
  \end{equation}Each of these recursions depends not just on other $T$ integrals, but also the $A$ and $B$ integrals. The $A$ and $B$ integrals can be split into line integrals over the sides $S_1, \ S_2, \ S_3$ of the triangle in the parameter space: $\int_{\partial \Omega_p} = \int_{S_1 \cup S_2 \cup S_3} = \int_{S_1}+\int_{S_2}+\int_{S_3}$. For $A_{m,n,q}$, this leads to 

  \begin{equation}
    \begin{split}
      A_{m,n,q}&= \int_{\partial \Omega_p} \alpha^m(\theta) \beta^n(\theta) \ R^{-q} \ \whvec{n}(\theta) \cdot \whvec e_1 \ d\theta \\
               &= \int_{S_1} \alpha^m(\theta) \beta^n(\theta) \ R^{-q} \ \whvec{n}(\theta) \cdot \whvec e_1 \ d\theta + \int_{S_2} \alpha^m(\theta) \beta^n(\theta) \ R^{-q} \ \whvec{n}(\theta) \cdot \whvec e_1 \ d\theta \\
      &\hspace{1em} + \int_{S_3} \alpha^m(\theta) \beta^n(\theta) \ R^{-q} \ \whvec{n}(\theta) \cdot \whvec e_1 \ d\theta\\
    \end{split}
    \label{eq:Amnq}
  \end{equation}These are integrals over line segments (see Figure \ref{fig:segment}), so $\whvec n(\theta) \cdot \whvec e_1\equiv \text{constant}$. For $B_{m,n,q}$, the only difference is $\whvec e_1$ is replaced by $\whvec e_2$. Each of these line segment integrals can be written as a linear combination of integrals that take one of the following forms,

    \begin{align}
  S^{\bell}_{m,-1}&:= \int_0^1 \theta^m R(\vec x_f, \vec y(\theta)) \ d \theta   \label{eq:Slm1} \\
  S^{\bell}_{m,1}&:= \int_0^1 \theta^m R^{-1}(\vec x_f, \vec y(\theta)) \ d \theta   \label{eq:Slp1}
    \end{align}The superscript $\bell$ denotes the direction parallel to the line segment over which we are integrating.

    Using this notation, we can write $A_{m,n,q}$ and $B_{m,n,q}$ as (see Appendix section \ref{ABderiv} for derivation)

    \begin{align}
      A_{m,n,q}&=S^{\whvec{e_2}}_{n,q} - \sum_{k=0}^{m+n}\binom{m+n}{k}(-1)^k S^{\whvec{d}}_{k,q} \label{eq:Amnq} \\
      B_{m,n,q}&=\begin{cases}
        \displaystyle - S^{\whvec{e_1}}_{m,q} + \sum_{k=0}^{m} \binom{m}{k}(-1)^{k}S^{\whvec{d}}_{k,q} & n=0 \\
                   \displaystyle \sum_{k=0}^{m+n} \binom{m+n}{k}(-1)^k S^{\whvec{d}}_{k,q} & n>0
                 \end{cases} \label{eq:Bmnq}
    \end{align}The $S$ integrals are evaluated by a recurrence formula derived in \cite{cortez18},

    \begin{equation}
    \label{eq:Srec}
    S^{\bell}_{m+1,q} = - \frac{1}{L^2(q-2)}\theta^mR^{-(q-2)} \Bigg \rvert_{\theta=0}^{\theta=1} + \frac{m}{L^2(q-2)}S^{\bell}_{m-1,q-2} - \frac{\vec x_0 \cdot \bell}{L}S_{m,q}^{\bell}, \ q \neq 2
  \end{equation}The base cases for the recursions with $q=\mp 1$ are,

  \begin{flalign}
    &S^{\bell}_{0,-1}=\int_0^1 R(\vec x_f, \vec y(\theta)) \ d\theta = \frac{1}{2L} \Bigg((\vec x(\theta) \cdot \bell)R(\vec x_f, \vec y(\theta)) \nonumber \\ 
             &{}  \hspace{3em}+ \Big(\vec x(\theta) \cdot \vec y(\theta))^2 - (\vec x(\theta) \cdot \bell)^2 +\epsilon^2\Big)\log \left [ \vec x(\theta) \cdot \bell + R(\vec x_f, \vec y(\theta))\right]\Bigg) \Bigg \rvert_{\theta=0}^{\theta=1} \label{eq:S0m1}   \\
    &S^{\bell}_{0,1} = \int_0^1 R^{-1}(\vec x_f, \vec y(\theta)) \ d\theta = \frac{1}{L} \ \text{arctanh}\left (\frac{\vec x(\theta) \cdot \bell}{R(\vec x_f, \vec y(\theta))} \right) \Bigg \rvert_{\theta=0}^{\theta=1} \label{eq:S0p1}
\end{flalign}
  \subsection{Base Cases $T_{0,0,3}$ and $T_{0,0,1}$}
The $T$ base cases that we need to use recursion formulas \eqref{eq:Tmp1}, \eqref{eq:Tnp1} are $T_{0,0,3}$ and $T_{0,0,1}$. In this section, we derive formulas for these integrals. 

\subsubsection{$T_{0,0,3}$}

\begin{figure}[ht]
    \centering
    \scalebox{0.7}{
    \def\svgwidth{\columnwidth}
    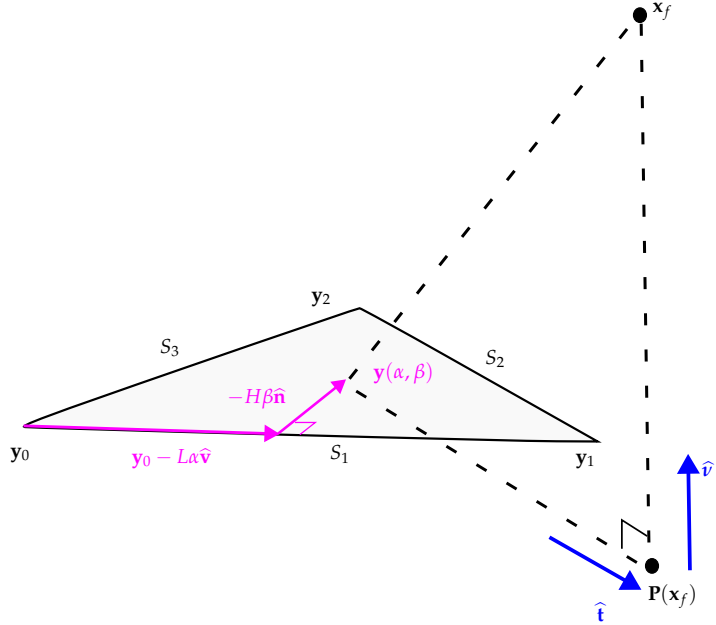
}
    \caption{A picture displaying the important geometric terms used in the derivation of $T_{0,0,3}$.}
    \label{fig:triPerspective}
\end{figure} Here, it will be convenient to use a parameterization different from the one shown in Figure \ref{fig:omega}. In general, for any point $\vec y$ lying in the triangular region, we can write its position in parametric form as, 

\[\vec y = \vec y(\alpha,\beta) = \vec y_0 - L\alpha \whvec v - H\beta \whvec n \]where $\whvec v$ is the unit vector in the direction of from the point $\vec y_1$ to $\vec y_0$, and $\whvec n$ is the outward pointing normal vector to the side $S_1$. Note that the physical coordinate system is now orthogonal. Continuing, for any field point $\vec x_f$ we can decompose the vector $\vec x_f - \vec y(\alpha,\beta)$ as  

\[ \vec x_f - \vec y(\alpha,\beta) = r \whvec t + z_0 \whvecg \nu \]where $\whvec t$ lies in the plane of the triangle and $\whvecg{\nu}$ is orthogonal to the plane of the triangle (see Figure \ref{fig:triPerspective}.) The squared regularized distance $R^2=\lvert  \vec x_f - \vec y(\alpha,\beta) \rvert^2 + \epsilon^2$ is

\[R^2=R^2(r)=r^2+z_0^2+\epsilon^2 = r^2+\gamma^2 \]where $\gamma^2=z_0^2+\epsilon^2$. The end goal is to convert the integral over the triangular region into a contour integral over its boundary. To this end, we will require the surface Laplacian operator $\Delta_{2D}$, which acting on a radially symmetric $C^2$ scalar function $f=f(r)$ is

\[\Delta_{2D}f:=\frac{1}{r}\pp{}{r}\left(r \pp{f}{r} \right)\]We note that for

\begin{equation}
  \label{eq:psi}
  \begin{split}
  \psi(r)&:=\frac{1}{\gamma} \text{log}\left ( R(r)+\gamma \right )\\
  \Delta_{2D} \psi &= \frac{1}{R^3} 
  \end{split}  
\end{equation}By Green's first identity, we have 

\[\int \int_{\Omega} \Delta_{2D} \left ( \psi(r) \right ) \ d \vec s =  \int \int_{\Omega} \frac{1}{R^3(r)} \ d \vec s = \oint_{\partial \Omega} \pp{\psi}{n} \ d\theta \]The rightmost integral above is the contour integral we will evaluate.  Since $\pp{\psi}{n} = \nabla \psi \cdot \whvec n$, we first need $\nabla \psi$. Note that all the derivatives mentioned here and in the following are with respect to the parameterized point $\vec y=\vec y(\alpha, \beta)$ lying on the triangle.

Using \eqref{eq:psi}, we have $\displaystyle \nabla \psi = \frac{1}{\gamma (R+\gamma)}\nabla R$ where $\displaystyle \nabla R = \frac{r \nabla r}{R}$ since $R^2= r^2 + \gamma^2$. The squared Euclidean distance in the plane of the triangle is

\[r^2 = \lvert \vec P(\vec x_f) - \vec y(\alpha,\beta) \rvert^2 = \lvert \vec P(\vec x_f) - \vec y_0 + L\alpha \whvec v + H\beta \whvec n \rvert^2\]where $\vec P(\vec x_f)=\vec x_f - z_0 \whvecg{\nu}$ is the projection of $\vec x_f$ onto the plane of the triangle. In terms of the parameterization,

\[r^2 = r^2(\alpha,\beta) = (\vec x_0 \cdot \whvec v + \alpha L)^2 + (\vec x_0 \cdot \whvec n + \beta H)^2  \]and its gradient $\nabla r$ (\textit{in physical space}) is 

\[\nabla r = \frac{1}{r} \begin{bmatrix}(\vec x_0 \cdot \whvec v + \alpha L)  \\ (\vec x_0 \cdot \whvec n + \beta H) \end{bmatrix}  \]Backsubstituting into $\nabla R$ and $\nabla \psi$, we have 

\[\nabla R = \frac{1}{R}\begin{bmatrix}(\vec x_0 \cdot \whvec v + \alpha L)\\ (\vec x_0 \cdot \whvec n + \beta H) \end{bmatrix} \]

\[\nabla \psi = \frac{1}{\gamma R (R+\gamma)}\begin{bmatrix}(\vec x_0 \cdot \whvec v + \alpha L)  \\ (\vec x_0 \cdot \whvec n + \beta H)  \end{bmatrix} \]Finally, we can use this to compute the part of $\displaystyle \oint_{\partial \Omega} \pp{\psi}{n} \ d\theta$ corresponding to the side $S_1$ from $\vec y_0$ to $\vec y_1$. In the new orthogonal coordinate system, $\whvec n$, the outward pointing normal to the triangle, is  $\begin{bmatrix}0 \\ -1 \end{bmatrix}$ and $\beta \equiv 0$. The integrand restricted to $S_1$ is

\begin{equation*}
  \begin{split}
    \nabla \psi \cdot \whvec n \Big \rvert_{S_1} &= \frac{1}{\gamma R (R+\gamma)}\begin{bmatrix}(\vec x_0 \cdot \whvec v + \alpha L) \\ (\vec x_0 \cdot \whvec n + \beta H)  \end{bmatrix} \cdot \begin{bmatrix} 0 \\ -1 \end{bmatrix} \Bigg \rvert_{\beta \equiv 0} \\
    &=-\frac{(\vec x_0 \cdot \whvec n)}{\gamma R(R + \gamma)}
  \end{split}
\end{equation*}Using $R(\alpha,\beta) \Big\rvert_{S_1}= \sqrt{(\vec x_0 \cdot \whvec v + \alpha L)^2+ (\vec x_0 \cdot \whvec n)^2 + \gamma^2}$, we have 

\begin{equation}
  \label{eq:S1int}
  \begin{split}
    &\int_{S_1} \nabla \psi \cdot \whvec n \ ds \\
    =& -\frac{\vec x_0 \cdot \whvec n}{\gamma L}\int_0^1 \frac{1}{\sqrt{(\alpha +P)^2+Q^2} \left(\sqrt{(\alpha +P)^2+Q^2}+\gamma/L \right)} \ d \alpha
  \end{split}
\end{equation}where $P=\frac{\vec x_0 \cdot \whvec v}{L}$ and $Q^2= (\frac{\vec x_0 \cdot \whvec n}{L})^2 + (\frac{\gamma}{L})^2$.

\subsubsection{Evaluation of integral \eqref{eq:S1int} }
Since $\gamma=\sqrt{z_0^2+\epsilon^2}$, the integrand of \eqref{eq:S1int} is always positive. In the case that $\vec x_0 \cdot \whvec n = 0$, \eqref{eq:S1int} evaluates to zero. On the other hand, if $\vec x_0 \cdot \whvec n \neq 0$, the integral evaluation depends on whether

\[-1 < P <0 \]The technical details of the following cases are covered in \ref{T003details} of the Appendix.

Suppose first that $P\geq 0$ or $P\leq -1$. Then it can be shown that 

\begin{align}
  \label{eq:T003sideeval1}
    &\int_0^1 \frac{1}{\sqrt{(\alpha +P)^2+Q^2} \left(\sqrt{(\alpha +P)^2+Q^2}+\gamma/L\right)} \ d \alpha \nonumber \\
 &=\begin{cases}
     \frac{2}{Q(1+\gamma/(LQ))} \sqrt{\frac{1+\gamma/(LQ)}{1-\gamma/(LQ)}} \tan^{-1}\left ( r\sqrt{\frac{1-\gamma/(LQ)}{1+\gamma/(LQ)}} \right ) \Bigg \rvert_{r=r_1}^{r=r_2}\text{if } P \geq 0 \\[1.5em]
    \frac{-2}{Q(1+\gamma/(LQ))} \sqrt{\frac{1+\gamma/(LQ)}{1-\gamma/(LQ)}} \tan^{-1}\left ( r\sqrt{\frac{1-\gamma/(LQ)}{1+\gamma/(LQ)}} \right ) \Bigg \rvert_{r=r_1}^{r=r_2}\text{if } P \leq -1 
   \end{cases}
  \end{align}where

  \begin{align*}
    r_1&= \tan \left (\frac{\sec^{-1}\left ( \sqrt{(P/Q)^2+1}\right )}{2}\right ) \\
    r_2 &= \tan \left ( \frac{\sec^{-1}\left ( \sqrt{(1+P)/Q)^2+1}\right )}{2}\right )
  \end{align*}On the other hand, if $-1<P<0$, then

\begin{equation}
  \label{eq:T003sideeval2}
    \begin{split}
      &\int_0^1 \frac{1}{\sqrt{(\alpha +P)^2+Q^2} \left(\sqrt{(\alpha +P)^2+Q^2}+\gamma/L\right)} \ d \alpha \\
      &=-\frac{1}{Q}\int_{r_1}^{0}\frac{1}{1+R/Q\frac{1-r^2}{1+r^2}} \frac{2}{1+r^2} \ dr + \frac{1}{Q}\int_{0}^{r_2}\frac{1}{1+\gamma/(LQ)\frac{1-r^2}{1+r^2}} \frac{2}{1+r^2} \ dr \ \\
      &= \frac{2}{Q(1+\gamma/(LQ))} \sqrt{\frac{1+\gamma/(LQ)}{1-\gamma/(LQ)}} \tan^{-1}\left ( r_1\sqrt{\frac{1-\gamma/(LQ)}{1+\gamma/(LQ)}} \right ) \\
      &+ \frac{2}{Q(1+\gamma/(LQ))} \sqrt{\frac{1+\gamma/(LQ)}{1-\gamma/(LQ)}} \tan^{-1}\left ( r_2\sqrt{\frac{1-\gamma/(LQ)}{1+\gamma/(LQ)}} \right ) \\
  \end{split}
\end{equation}To compute $T_{0,0,3}$, the process is repeated for each side. This can be summarized neatly as an algorithm (Algorithm \ref{alg:T003}.) Here, suppose that each side $S_i, \ i=1,2,3$ is a data structure with two components containing the vertices of the side ordered in the way that the triangle is traversed. Using our notation, we have

\begin{align*}
  S_1: S_1(1)=\vec y_0, \ S_1(2)=\vec y_1 \\
  S_2: S_2(1)=\vec y_1, \ S_2(2)=\vec y_2 \\
  S_3: S_3(1)=\vec y_2, \ S_3(2)=\vec y_0 
\end{align*}

\begin{algorithm}[t!]
  \caption{$T_{0,0,3}$}\label{alg:T003}
  \begin{algorithmic}
    \Require{Field point $\vec x_f$ and triangle with sides $S_1,S_2,S_3$}
    \State $T_{0,0,3} \gets 0$ \Comment{Initialize to zero}
    \For{$i=1:3$}
    \State $(\vec y_a, \vec y_b) \gets (S_i(1),S_i(2))$
    \State $\vec x_0 \gets \vec x_f - \vec y_a$
    \State $\whvec n \gets$ outward unit normal vector to $S_i$
    \If{$(\vec x_0 \cdot \whvec n) \neq 0$} 
    \State $L \gets \lvert \vec y_a -\vec y_b \rvert $ 
    \State $\whvec v \gets \frac{\vec y_a - \vec y_b}{L}$
    \State $H \gets$ projection of $S_j$ onto $\whvec n$, $j\neq i$
    \State $P \gets \frac{\vec x_0 \cdot \whvec v}{L}$
    \State $Q^2 \gets (\frac{\vec x_0 \cdot \whvec n}{L})^2
    + (\frac{\gamma}{L})^2$
    \If{$-1<P<0$}
    \State $\text{out}_i \gets - \frac{\vec x_0 \cdot \whvec n}
    {\gamma L}\times$ \eqref{eq:T003sideeval2}
    \Else
    \State $\text{out}_i \gets - \frac{\vec x_0 \cdot \whvec n}
    {\gamma L} \times $ \eqref{eq:T003sideeval1}
    \EndIf
    \State $T_{0,0,3} \gets T_{0,0,3}+\text{out}_i$ 
    \EndIf
    \EndFor
    \State $T_{0,0,3} \gets T_{0,0,3} / (LH)$ \Comment{for consistency with formula \eqref{eq:velotri}}
    \State \Return $T_{0,0,3}$
  \end{algorithmic}
\end{algorithm}

\subsubsection{$T_{0,0,1}$}
The same parameterization and notation described in the previous section will be used in the computation of the integral $\displaystyle T_{0,0,1}=\int \int_{\Omega}\frac{1}{R} \ d\vec s$. We note that the surface Laplacian of $R(r)$ is,

\[\Delta_{2D}(R) = \frac{1}{R}+\frac{\gamma^2}{R^3}\]Applying Green's first identity, we have 

\[\int \int_{\Omega} \Delta_{2D} \left ( R(r) \right ) \ d \vec s =  \int \int_{\Omega} \frac{1}{R(r)} + \frac{\epsilon^2}{R^3(r)} \ d \vec s = \oint_{\partial \Omega} \pp{R}{n} \ d\theta \]Rearranging, this is equivalent to 

\[\int \int_{\Omega} \frac{1}{R} \ d\vec s = \oint_{\partial \Omega}\pp{R}{n}  \ d\theta -\epsilon^2 \int \int_{\Omega}\frac{1}{R^3}\]where $\pp{R}{n}$ is the derivative of $R$ in the direction normal to the boundary of the triangle and the last integral on  the right hand side is $T_{0,0,3}$. The contour integral can be evaluated side by side. We use the first side, $S_1$ as an example as before. From the previous section, we have 

\[ \nabla R = \frac{1}{R}\begin{bmatrix}(\vec x_0 \cdot \whvec v + \alpha L)\\ (\vec x_0 \cdot \whvec n + \beta H) \end{bmatrix} \]We also have that $\whvec n =\begin{bmatrix}0 \\ -1 \end{bmatrix}$ and $\beta \equiv 0$ along $S_1$. This simplifies $\pp{R}{n} \Big \rvert_{S_1}$ as

\begin{align}
  \pp{R}{n} \Big \rvert_{S_1} &= \nabla R \cdot \whvec n \Big \rvert_{S_1}= \nabla R \Big \rvert_{\beta=0} \cdot \begin{bmatrix}0 \\ -1 \end{bmatrix} \nonumber \\
  &= -\frac{\vec x_0 \cdot \whvec n}{\sqrt{(\vec x_0 \cdot \whvec v + \alpha L)^2+ (\vec x_0 \cdot \whvec n)^2 +\gamma^2}}
\end{align}We then have

\begin{equation}
  \label{eq:T001side}
  \begin{split}
  \int_{S_1} \pp{R}{n} \ d \vec s &= -(\vec x_0 \cdot \whvec n) \int_0^1 \frac{d\alpha}{\sqrt{(\vec x_0 \cdot \whvec v + \alpha L)^2 + (\vec x_0 \cdot \whvec n)^2+ \gamma^2}} \\
  &=  - \frac{\vec x_0 \cdot \whvec n }{L} \text{arctanh}\left (\frac{\vec x(\alpha) \cdot \whvec v}{R(\vec x_f, \vec y(\theta))} \right) \Bigg \rvert_{\alpha=0}^{\alpha=1}  
  \end{split}  
\end{equation}where $\vec x(\alpha) = \vec x_f - \vec y(\alpha)$. This formula is not new: it was used before in \eqref{eq:S0p1}. Having computed $T_{0,0,3}$, we can now compute $T_{0,0,1}$ using Algorithm \ref{alg:T001}.

\begin{algorithm}
  \caption{$T_{0,0,1}$}\label{alg:T001}
  \begin{algorithmic}
    \Require{Field point $\vec x_f$, triangle with sides $S_1,S_2,S_3$, and $T_{0,0,3}$}
    \State $T_{0,0,1} \gets 0$ \Comment{Initialize to zero}
    \For{$i=1:3$}
    \State $(\vec y_a, \vec y_b) \gets (S_i(1),S_i(2))$
    \State $\vec x_0 \gets \vec x_f - \vec y_a$
    \State $\whvec n \gets$ outward unit normal vector to $S_i$
    \If{$(\vec x_0 \cdot \whvec n) \neq 0$} 
    \State $L \gets \lvert \vec y_a -\vec y_b \rvert $ 
    \State $\whvec v \gets \frac{\vec y_a - \vec y_b}{L}$
    \State $H \gets$ projection of $S_j$ onto $\whvec n$, $j\neq i$
    \State $T_{0,0,1} \gets T_{0,0,1} +$ \eqref{eq:T001side}
    \EndIf
    \EndFor
    \State $T_{0,0,1} \gets T_{0,0,1} - \epsilon^2 T_{0,0,3}$
    \State $T_{0,0,1} \gets T_{0,0,1}/ (LH)$ \Comment{for consistency with formula \eqref{eq:velotri}} 
    \State \Return $T_{0,0,1}$
  \end{algorithmic}
\end{algorithm}

\subsection{Summary}
For each triangle and evaluation point, there are eight base cases that must be established before using the recurrence formulas \eqref{eq:Tmp1} and \eqref{eq:Tnp1}. These are the following,

\[S^{\whvec{e_1}}_{0,-1}, \ S^{\whvec{e_1}}_{0,1}, \ S^{\whvec{e_2}}_{0,-1}, \ S^{\whvec{e_2}}_{0,1}, \ S^{\whvec{d}}_{0,-1}, \ S^{\whvec{d}}_{0,1}, \ T_{0,0,3}, \ T_{0,0,1} \]With these base cases in hand, we use \eqref{eq:Srec} to compute the $A_{m,n,q}$ and $B_{m,n,q}$, followed by \eqref{eq:Tmp1}, \eqref{eq:Tnp1} to compute the $T_{m,n,q}$ in formula \eqref{eq:velotri}.

\section{Examples}
In the following examples, we will measure the errors using the Euclidean vector norm when measuring a single quantity, and the $\ell_2$ norm for the average error of several measurements. Mathematically, for any point $\vec x$, we define the error $e(\vec x)$ in the numerically computed velocity field $\vec u_{\text{numerical}}$ as 

\[e (\vec x) \equiv \lvert \vec u_{\text{exact}} - \vec u_{\text{numerical}} \rvert \]In the case of several points, $\vec x_1, \vec x_2, \dots, \vec x_n$, we have the associated errors $e(\vec x_1), e(\vec x_2),$ \\ $\dots e(\vec x_n)$, which we will write concisely as $e_1, e_2, \dots, e_n$. The $\ell_2$ norm of these errors is defined,

\begin{align}
  \|{e}\|_{2} &= \sqrt{\frac{1}{n}\sum_{k=1}^n e^2_k}
\end{align}

\subsection{Uniformly translating/rotating sphere}
The velocity field due to a uniformly translating and/or rotating sphere in an infinite fluid is a well known result \cite{landau}. As a test problem, it can be seen from two vantage points:

\begin{enumerate}
\item Given the traction field on the sphere, triangulate the sphere and use the forces at the vertices to approximate the velocity field at several points. Then compare this velocity field to the theoretical velocity field.
\item Given the velocity field on the sphere, triangulate the sphere and solve a linear system to compute the forces at the vertices. Then compare the net force or torque to the theoretical results.  
\end{enumerate}We will call the first of these the \textit{forward problem} and the second the \textit{resistance problem.} In the following, we will discuss results from each approach. But before this discussion, we describe a particular triangulation of the sphere that we utilize frequently.

\subsubsection{Delaunay triangulation}
To triangulate the sphere, we take points on a regular icosahdron with vertices that coincide with the unit sphere. The triangular faces are subdivided with the number of subdivisions specified by a positive integer factor $f$. For a given $f$, the number of triangular faces generated by this procedure is $20f^2$. The new vertices generated are projected onto the unit sphere and the triangulation is generated by the convex hull of the projected points. This is equivalent to a Delaunay triangulation and the algorithm described is implemented in the MATLAB library \textit{SPHERE\_DELAUNAY} published by John Burkardt \cite{burkardt}.

\subsubsection{The forward problem}
\noindent \textit{The translating sphere} \\
For a sphere of radius $a$ uniformly translating in an infinite fluid with velocity $\vec U$, the hydrodynamic traction $\vec f$ at every point on the sphere is

\begin{equation}
  \label{eq:sphereTraction}
  \vec f = -\frac {3 \mu}{2 a} \vec U 
\end{equation}The velocity field at any point in the fluid, $\vec u$, is axisymmetric about the translational direction. Using the center of the sphere as the origin, $\vec u(\vec x)$ is written ($r=\lvert \vec x \rvert)$

\begin{equation}
  \label{eq:uSphere}
  \begin{split}
    \vec u (\vec x) &= \frac{a}{4r}\left (3 + \frac{a^2}{r^2} \right ) \vec U + \frac{3a (\vec x \cdot \vec U)}{4r^2}  \left (1 - \frac{a^2}{r^2}\right )\frac{\vec x}{r}
  \end{split}
\end{equation}We test the method by imposing the traction \eqref{eq:sphereTraction} on the vertices of the triangle and comparing the computed velocity field to the theoretical result \eqref{eq:uSphere}. In the following, we take $a=1$, $\vec U=(1,0,0)$, and $\mu=1$.

In Table \ref{Tab:forwardSurfaces}, we tabulate the $\ell_2$ errors for several discretizations and regularizations. We also record the number of triangles (``no. triangles'') and the number of degrees of freedom (DOF = 3 $\times$ (number of discretization nodes).) The discretization $h$ is the square root of twice the average area of a triangle ($h=\sqrt{\text{mean}\left(BH\right)}$).

\begin{table}[ht]
  \begin{subtable}[h]{0.9\textwidth} 
\centering
  \begin{tabular}{l|l|l?c|c|c|c|c|c}%
    \hline\hline & & & \multicolumn{6}{c}{\bfseries Regularization ($\epsilon$)}\\
     \bfseries no. tris & \bfseries DOF & \bfseries Disc. ($h$) & \bfseries $10^{-1}$ & $10^{-2}$ & $10^{-4}$ & $10^{-6}$ & $10^{-7}$ & $10^{-8}$
                                                                                                                                                           \begin{filecontents*}{tables/forwardSurfacesCollocation.csv}
numFaces,DOF,discretization,epsone,epstwo,epsthree,epsfour,epsfive,epssix
180,276,0.3675,0.0532,0.0285,0.0257,0.0257,0.0257,0.0257
720,1086,0.1860,0.0351,0.0093,0.0064,0.0064,0.0064,0.0064
2880,4326,0.0933,0.0307,0.0046,0.0016,0.0016,0.0016,0.0016
11520,17286,0.0467,0.0296,0.0034,0.0004,0.0004,0.0004,0.0004
\end{filecontents*}
\csvreader[head to column names]{tables/forwardSurfacesCollocation.csv}{}
                                                                                                                                                           {\\\hline  \numFaces & \DOF & \discretization & \epsone & \epstwo  & \epsthree & \epsfour & \epsfive & \epssix}
  \end{tabular}
  \caption{Stokeslet surfaces $\ell_2$ error in velocity at triangle vertices in the forward problem for translating sphere.}
  \label{Tab:forwardSurfaces}
  \end{subtable}\\
\begin{subtable}[h]{0.9\textwidth}
\centering
  \begin{tabular}{l|l|l?c|c|c|c|c|c}%
    \hline\hline & & & \multicolumn{6}{c}{\bfseries Regularization ($\epsilon$)}\\
     \bfseries no. tris & \bfseries DOF & \bfseries Disc. ($h$) & \bfseries $0.25$ & $0.1$ & $0.05$ & $0.01$ & $0.005$ & $0.001$
                                                                                                                                                           \begin{filecontents*}{tables/forwardMRSCollocation.csv}
numFaces,DOF,discretization,epsone,epstwo,epsthree,epsfour,epsfive,epssix
180,276,0.3675,0.0927,0.0271,0.1736,1.4377,3.0192,15.6719
720,1086,0.1860,0.0764,0.0310,0.0203,0.3456,0.7570,4.0484
2880,4326,0.0933,0.0711,0.0309,0.0142,0.0730,0.1767,1.0079
11520,17286,0.0467,0.0698,0.0297,0.0153,0.0113,0.0371,0.2454
\end{filecontents*}
\csvreader[head to column names]{tables/forwardMRSCollocation.csv}{}
                                                                                                                                                           {\\\hline \numFaces & \DOF & \discretization & \epsone & \epstwo  & \epsthree & \epsfour & \epsfive & \epssix}
  \end{tabular}
  \caption{MRS $\ell_2$ error in velocity at triangle vertices in the forward problem for translating sphere.}
  \label{Tab:forwardMRS}
\end{subtable}
\caption{Errors for the translating sphere problem}
\label{Tab:trans}
\end{table}

\begin{figure}[h!]
  \centering
  \begin{subfigure}{0.49\textwidth}
    {\includegraphics[width=\textwidth]{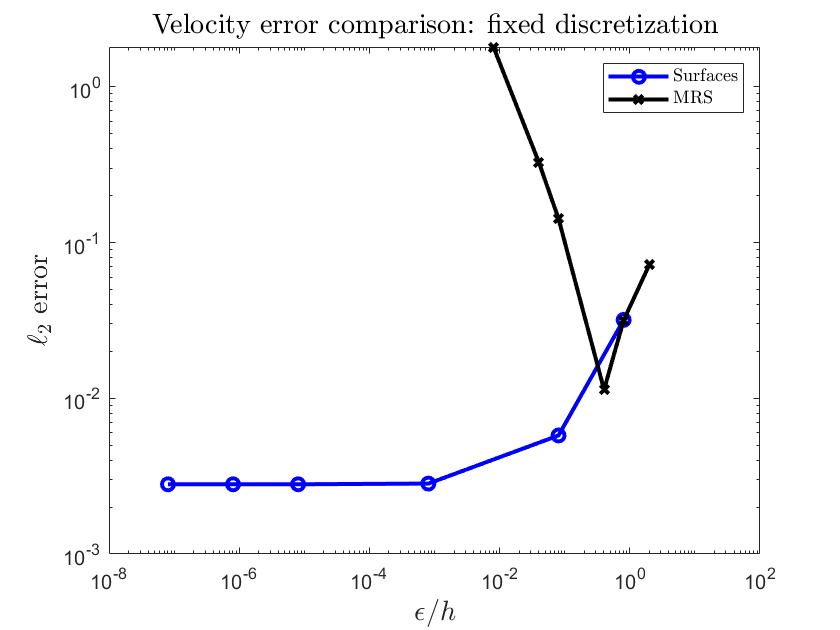}}
    \caption{}
    \label{fig:tsfixh}
  \end{subfigure}
  \begin{subfigure}{0.49\textwidth}
    {\includegraphics[width=\textwidth]{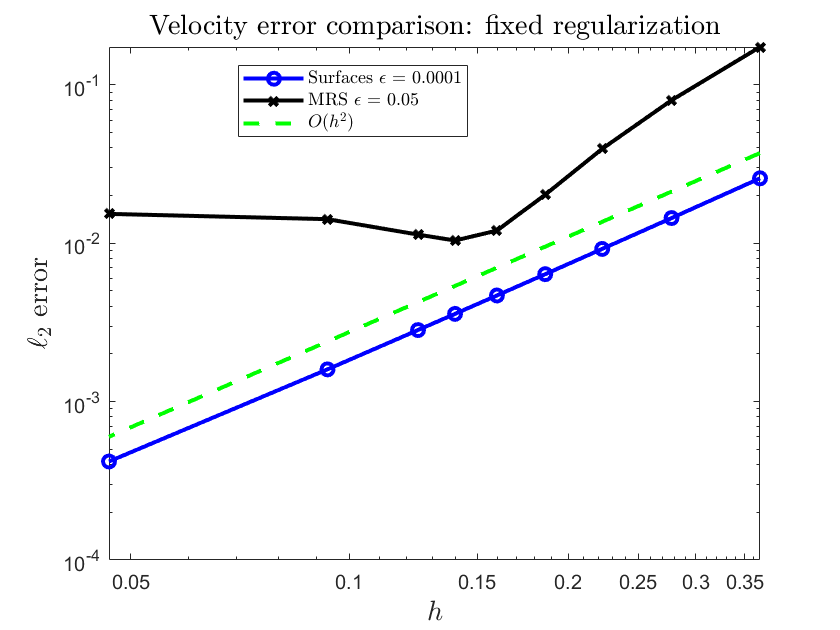}}
    \caption{}
    \label{fig:tsfixeps}
  \end{subfigure}
  \caption{A comparison of Stokeslet surfaces and MRS $\ell_2$ errors for the forward problem of computing the velocity on a sphere using the hydrodynamic traction $\vec f = - \frac{3 \mu}{2a}(1,0,0).$ (a) Error as a function of $\epsilon/h$. The spatial discretization is fixed for both methods as $h=0.1243$. (b) Error as a function of $h$. For the MRS, $\epsilon$ was set to $0.05$ while for the Stokeslet surface method it was set to $10^{-4}$.}
    \label{fig:ts}
\end{figure}

As a comparison, the results using the method of regularized Stokeslets (MRS) are tabulated in Table \ref{Tab:forwardMRS}. For MRS, we use the same spatial discretizations as those used with the surfaces, but relatively larger regularizations. This is due to the fact that the quadrature error is $O(h^2/\epsilon^3)$ \cite{cortez2005} so the error grows when $\epsilon \ll h$. This behavior can be seen in Figure \ref{fig:tsfixh}: the $\ell_2$ error decreases as long as $\epsilon/h > 10^{-1}$ but grows beyond this threshold. On the other hand, the error for the Stokeslet surfaces decreases until $\epsilon/h \approx 10^{-3}$ and remains the same even as $\epsilon$ gets smaller relative to $h$.

If $\epsilon$ is fixed but the spatial discretization $h$ is varied, we observe second order convergence using the Stokeslet surfaces. This is apparent in the log-log plot in Figure \ref{fig:tsfixeps} with data falling parallel to the dashed green line indicating second order convergence. The regularization used for the Stokeslet surfaces in this figure is $\epsilon = 10^{-4}$. In contrast, the results from MRS using $\epsilon = 0.05$ appear to converge at an approximately second order rate until $h=0.1596$. For $h$ smaller than this threshold, the error slightly increases but then flattens out for the finer discretizations, indicating the dominance of the regularization error.

\begin{figure}[htbp!]
  \centering
  \begin{subfigure}{0.4\textwidth}
    \centering
    \includegraphics[width=\textwidth]{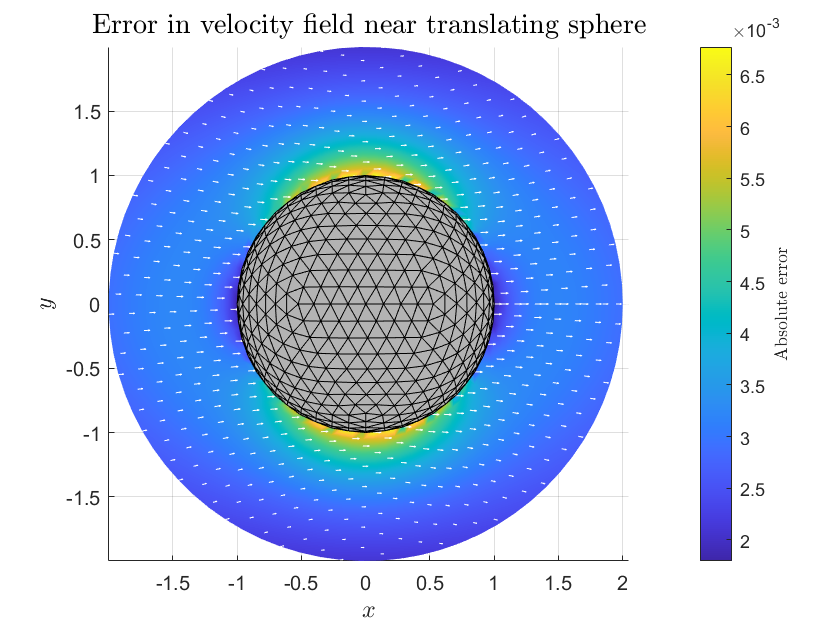}
    \caption{Absolute error}
  \end{subfigure}
  \begin{subfigure}{0.4\textwidth}
    \centering
  \includegraphics[width=\textwidth]{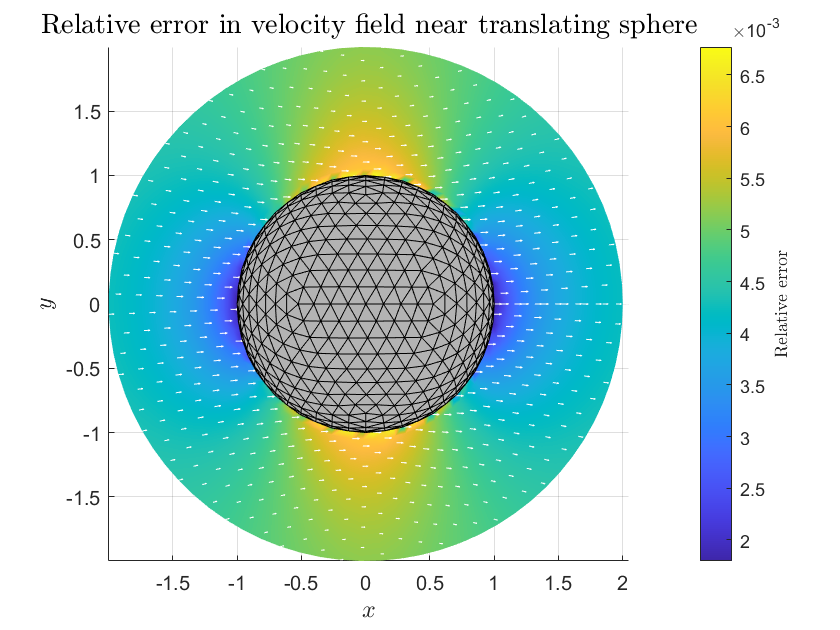}
  \caption{Relative error}
\end{subfigure}
\caption[Visualization of velocity field and errors near translating sphere]{A visualization of the error in the velocity field for the unit sphere (spatial discretization $h=0.1398$, regularization $\epsilon=10^{-4}$) with imposed traction $\vec f=-\frac{3 \mu}{2a}(1,0,0)$. The traction induces rigid translation of the sphere with unit velocity in the x direction. The errors are shown in the equatorial plane ($z=0$). (a) illustrates that the largest errors (by absolute magnitude) are concentrated near the sphere on the sides of the sphere orthogonal to the translational direction. (b) shows the largest relative errors are about 0.6\% in the yellow band  orthogonal to the translational direction. White arrows show the fluid velocity in the lab frame.}
\label{fig:translationColor}
\end{figure}

  A visualization of the errors in the vicinity of the sphere is shown in Figure \ref{fig:translationColor} for the particular case of $h=0.1398, \ \epsilon=10^{-4}$. The errors are color coded with ``warm'' colors corresponding to relatively large error magnitudes and ``cool'' colors corresponding to relatively small error magnitudes. In Figure \ref{fig:translationColor} (a), the absolute error is used. External to the sphere, the errors are largest on the ``sides'' of the sphere that are orthogonal or nearly orthogonal to the direction of motion. In (b), we visualize the relative error in velocity where

\[\text{relative error} = \frac{\lvert \vec U_{\text{approx}} - \vec U_{\text{theory}} \rvert}{\lvert \vec U_{\text{theory}}\rvert} \]The largest relative errors ($\approx 1\%$) are concentrated in the yellow band orthogonal to the direction of motion.

\noindent \textit{The rotating sphere} \\
In the case of a sphere uniformly rotating with angular velocity $\vecg \Omega$, the hydrodynamic traction $\vec f$ for a point $\vec x$ on the sphere is

\begin{equation}
  \vec f(\vec x) = \frac{3 \mu}{a} \  \vec x \times \vecg \Omega 
  \label{eq:tractionRotation}
\end{equation}The resultant velocity field, $\vec u$, is

\begin{equation}
  \vec u = \frac{a^3}{r^2}\sin \phi \ \vecg \Omega \times \vec x 
  \label{eq:rotationVelocity}
\end{equation}where $\phi$ is the angle between the ray connecting the center of the sphere to $\vec x$ and the axis of rotation. In the following discussion, we will fix $a=1$, $\mu=1$, and $\vecg \Omega = (0,0,1)$ so that we have a unit sphere rotating counterclockwise with unit angular speed about the z-axis.

Since the traction is no longer constant on the sphere, we expect there to be an advantage in the linearization of the forces that occurs automatically with the method of Stokeslet surfaces. To test this, instead of comparing the method to MRS like in the previous section, we compare to a version of the Stokeslet surfaces where the force distribution is constant. In particular, for a triangle with vertices $\vec y_0, \vec y_1, \vec y_2$ on the sphere, we prescribe the hydrodynamic traction $\vec f_{\text{cons}}$

\[\vec f_{\text{cons}} = \frac{1}{3} \left (\vec f(\vec y_0) + \vec f(\vec y_1) + \vec f(\vec y_2)\right )  \]where $\vec f(\vec y_j)$ comes from the traction formula \eqref{eq:tractionRotation}. The force density $\vec f_{\text{cons}}$ over each triangle is then just a constant function that averages the exact traction at the three vertices.

Like in the translating sphere problem, we compute the $\ell_2$ error from the computed velocity and the theoretical velocity $\vecg \Omega \times \vec x$ for all of the triangular vertices. The results for the piecewise linear forces are shown in Table \ref{Tab:forwardLinearSurfacesTorque} and the results for the piecewise constant forces are shown in Table \ref{Tab:forwardConstantSurfacesTorque}.

\begin{table}[ht]
  \begin{subtable}[h]{0.9\textwidth}
    \centering
    \begin{tabular}{l|l|l?c|c|c|c|c|c}%
      \hline\hline & & & \multicolumn{6}{c}{\bfseries Regularization ($\epsilon$)}\\
      \bfseries no. tris & \bfseries DOF & \bfseries Disc. ($h$) & \bfseries $10^{-1}$ & $10^{-2}$ & $10^{-4}$ & $10^{-6}$ & $10^{-7}$ & $10^{-8}$
                                                                                                                                         \begin{filecontents*}{tables/forwardSurfacesTorque.csv}
numFaces,DOF,discretization,epsone,epstwo,epsthree,epsfour,epsfive,epssix
180,276,0.3675,0.0958,0.0438,0.0381,0.0380,0.0380,0.0380
720,1086,0.1860,0.0697,0.0154,0.0094,0.0094,0.0094,0.0094
2880,4326,0.0933,0.0632,0.0084,0.0024,0.0023,0.0023,0.0023
11520,17286,0.0467,0.0616,0.0067,0.0006,0.0006,0.0006,0.0006
\end{filecontents*}
\csvreader[head to column names]{tables/forwardSurfacesTorque.csv}{}
                                                                                                                                         {\\\hline \numFaces & \DOF & \discretization & \epsone & \epstwo  & \epsthree & \epsfour & \epsfive & \epssix}
    \end{tabular}
    \caption{Piecewise linear Stokeslet surfaces $\ell_2$ error in velocity at triangle vertices in the forward problem for rotating sphere.}
    \label{Tab:forwardLinearSurfacesTorque}
  \end{subtable}
  \begin{subtable}[h]{0.9\textwidth}
    \centering 
    \begin{tabular}{l|l|l?c|c|c|c|c|c}%
      \hline\hline & & & \multicolumn{6}{c}{\bfseries Regularization ($\epsilon$)}\\
      \bfseries no. tris & \bfseries DOF & \bfseries Disc. ($h$) & \bfseries $10^{-1}$ & $10^{-2}$ & $10^{-4}$ & $10^{-6}$ & $10^{-7}$ & $10^{-8}$
                                                                                                                                         \begin{filecontents*}{tables/forwardConstantSurfacesTorque.csv}
numFaces,DOF,discretization,epsone,epstwo,epsthree,epsfour,epsfive,epssix
180,540,0.3675,0.1062,0.0551,0.0495,0.0495,0.0495,0.0495
720,2160,0.1860,0.0722,0.0182,0.0123,0.0122,0.0122,0.0122
2880,8640,0.0933,0.0638,0.0091,0.0031,0.0030,0.0030,0.0030
11520,34560,0.0467,0.0617,0.0069,0.0008,0.0007,0.0007,0.0007
\end{filecontents*}
\csvreader[head to column names]{tables/forwardConstantSurfacesTorque.csv}{}
                                                                                                                                         {\\\hline \numFaces & \DOF & \discretization & \epsone & \epstwo  & \epsthree & \epsfour & \epsfive & \epssix}
    \end{tabular}
    \caption{Piecewise constant Stokeslet surfaces $\ell_2$ error in velocity at triangle vertices in the forward problem for rotating sphere.}
    \label{Tab:forwardConstantSurfacesTorque}
  \end{subtable}
  \caption{Errors for the rotating sphere problem.}
  \label{Tab:rot}
\end{table}

\begin{figure}[h!]
  \centering
  \begin{subfigure}{0.49\textwidth}
    \includegraphics[width=\textwidth]{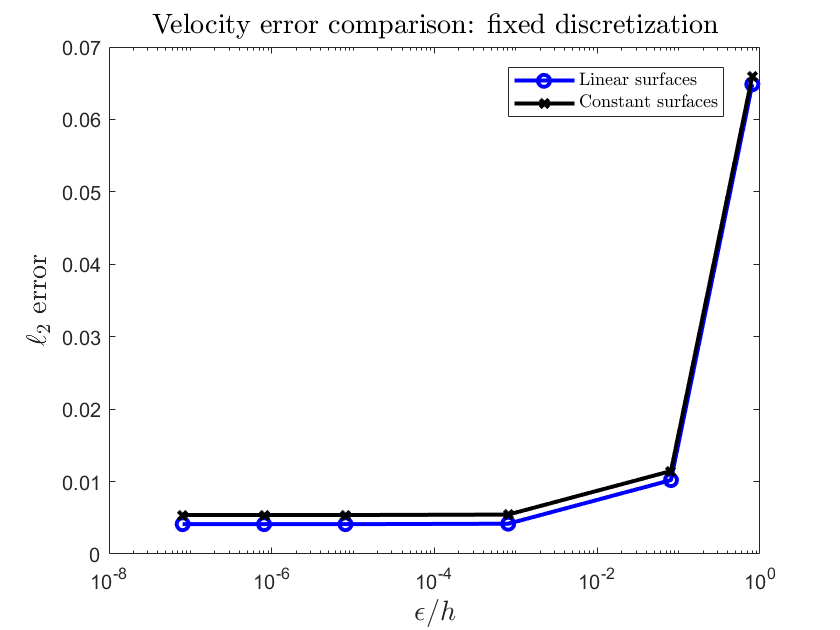}
    \caption{}
    \label{fig:lvsca}
  \end{subfigure}
  \begin{subfigure}{0.49\textwidth}
    \includegraphics[width=\textwidth]{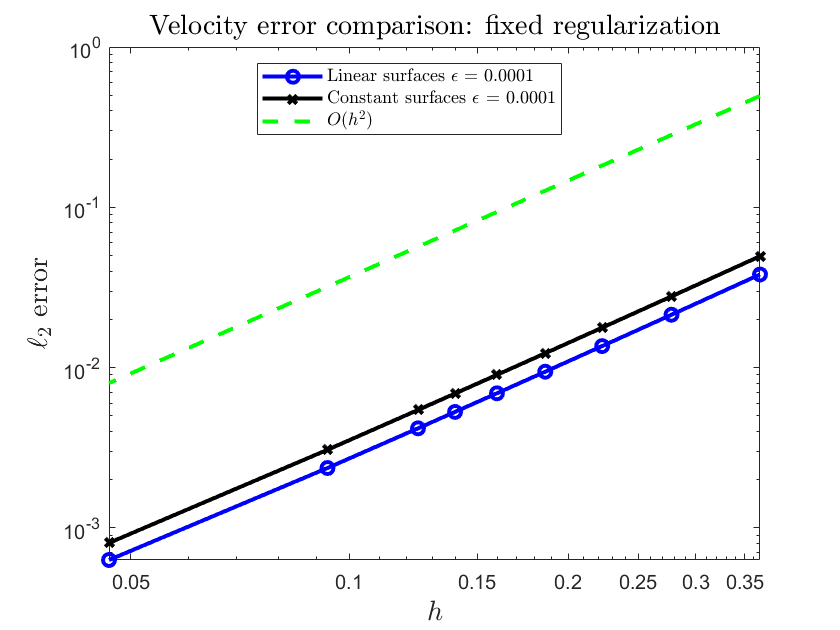}
    \caption{}
    \label{fig:lvscb}
  \end{subfigure}
  \caption{A comparison of linear Stokeslet surfaces and constant Stokeslet surfaces for the forward problem of computing the velocity on the unit sphere uniformly rotating with angular velocity $\vecg \Omega=(0,0,1)$. (a) Error as a function of $\epsilon/h$. The spatial discretization is fixed for both methods as $h=0.1243$. (b) Error as a function of $h$. The regularization, $\epsilon$, is a fixed value for each method while the discretization, $h$, is varied.}
  \label{fig:linearVsConstant}
\end{figure}

The tables show improvement in the evaluation of the velocity field from linearizing the force density compared to using a constant force density, with the largest differences in performance appearing for the coarser discretizations. Figure \ref{fig:lvsca} shows that by fixing the discretization and varying $\epsilon$,  the error curves for each method flatten out when $\epsilon/h \leq 10^{-3}$, indicating the dominance of the discretization error over the regularization error for values of $\epsilon/h$ less than this threshold. On the other hand, when the regularization is fixed, both methods decay at a second order rate in the spatial discretization $h$ (Figure \ref{fig:lvscb}.)

\begin{figure}[h!]  \centering
  \begin{subfigure}{0.48\textwidth}
{\includegraphics[width=\textwidth]{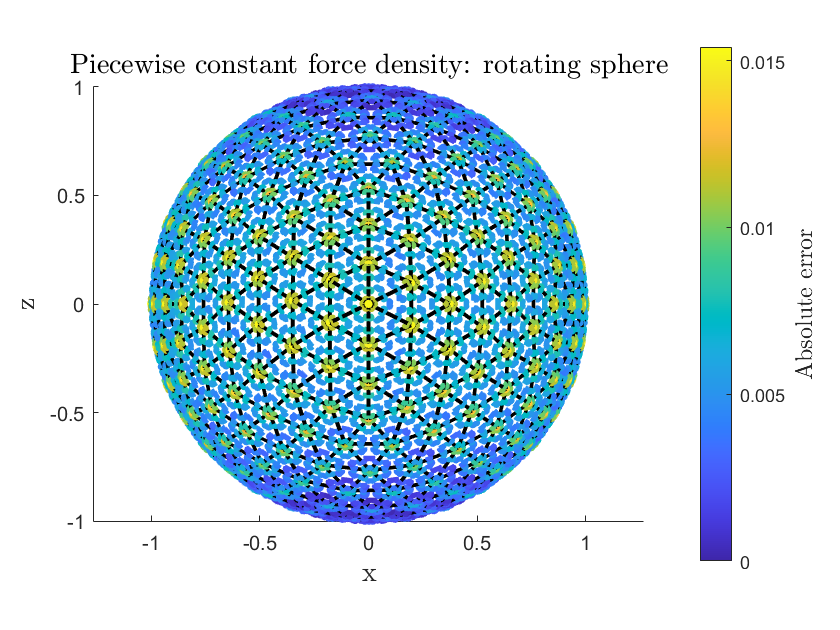}}
    \caption{}
    \label{fig:consside}
  \end{subfigure}
  \begin{subfigure}{0.48\textwidth}
{\includegraphics[width=\textwidth]{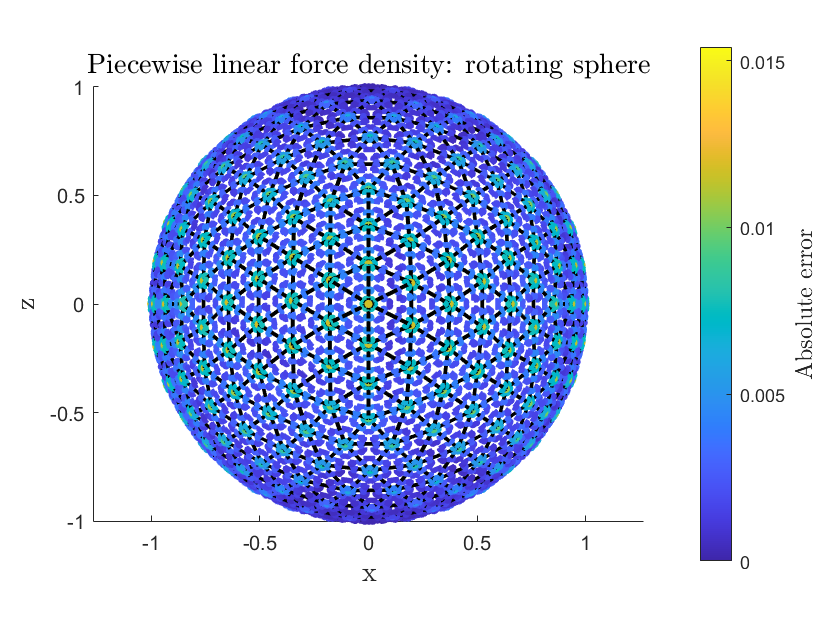}}
    \caption{}
    \label{fig:linside}
  \end{subfigure}\\
  \begin{subfigure}{0.48\textwidth}
{\includegraphics[width=\textwidth]{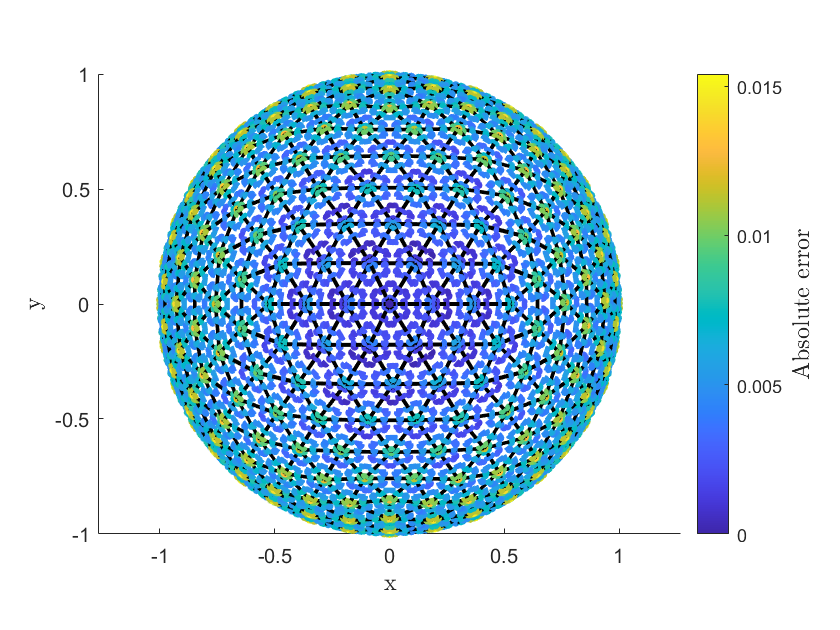}}
    \caption{}
    \label{fig:constop}
  \end{subfigure}
  \begin{subfigure}{0.48\textwidth}
{\includegraphics[width=\textwidth]{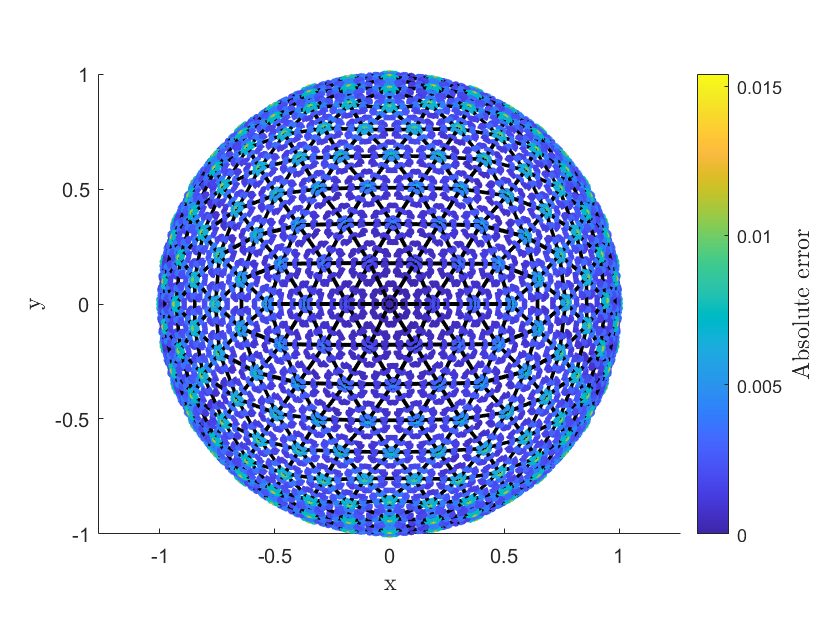}}
    \caption{}
    \label{fig:lintop}
  \end{subfigure}
  \caption{Performance comparison of piecewise linear forces and
piecewise constant forces for the problem of a unit sphere rotating with unit angular velocity about the positive z-axis. Each sphere is triangulated by 720 triangles ($h\approx 0.1860$) and the regularization is fixed at $\epsilon=10^{-4}$. The evaluation points lie on the triangles and the velocity at each point is compared to the theoretical velocity that ensures rigid rotation of the body. (a/c) show the side/top view for the rotating sphere using a piecewise constant force density. (b/d) show the side/top view for the rotating sphere using a piecewise linear force density.}
    \label{fig:rotSphereComp}
\end{figure}

Figure \ref{fig:rotSphereComp} illustrates the differences in performance on the surface of the rotating sphere for a specific case ($h\approx 0.1860$, $\epsilon=10^{-4}$). The evaluation points are the triangle vertices and several points in the interior of each triangle. Each evaluation point is represented as a dot and the errors are color coded by the absolute error. The errors are noticeably larger near the triangle vertices when using the piecewise constant force density. On the other hand, the error magnitudes are smaller and more uniform when using the linear force density. This difference is most apparent near the equator (the side views \ref{fig:consside}, \ref{fig:linside}) where the traction \eqref{eq:tractionRotation} is largest.

\subsubsection{The resistance problem}
In this section, we take the second vantage point explained earlier. Instead of using a given traction on the triangulated sphere and directly computing the resultant velocity, we impose the velocity at boundary points and solve a $3N \times 3N$ linear system for the force density at the $N$ triangle vertices. We can then compute the hydrodynamic drag (translating sphere) or moment (rotating sphere) using \eqref{eq:netForce} or \eqref{eq:netTorque} and compare with the theoretical results. \\[0.1em]

\noindent \textit{Drag on a translating sphere}\\
The theoretical value for the drag on a translating sphere, first computed by Stokes, is

\[\vec F_{\text{drag}} = - 6 \pi \mu a \vec U \]where $a$ is the radius of the sphere and $\vec U$ is the translational velocity.

\begin{figure}[h!]
  \centering
  \includegraphics[width=0.6\textwidth]{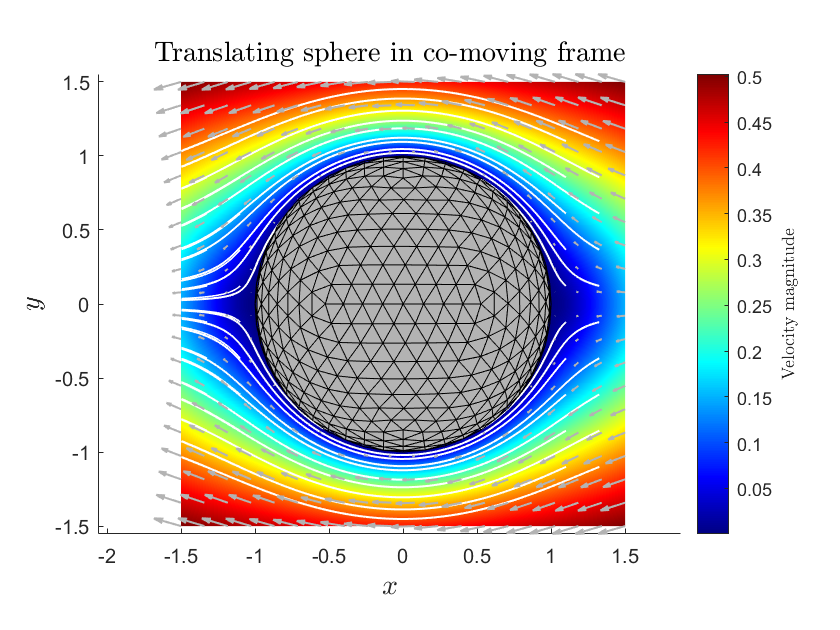}
  \caption[Flow around translating sphere in a co-moving frame.]{Flow around translating sphere in a co-moving frame. The white lines shown are streamlines and the magnitude of the velocity is shown in color.}
  \label{fig:sphereFlowAround}
\end{figure}

We use the unit sphere ($a=1$) and fix the viscosity $\mu=1$. At all of the triangle vertices, we impose the velocity $\vec U=(1,0,0)$. The chosen parameters give a theoretical drag of $\vec F_{\text{drag}}=(-6 \pi,0,0)$. The numerical drag, $\widetilde{\vec F}_{\text{drag}}=(\tilde{f}_1,\tilde{f}_2,\tilde{f}_3)$, is computed by summing the net force from every triangle with equation \eqref{eq:netForce}. A visualization of the flow around the sphere in a co-moving frame from a sample computation is shown in Figure \ref{fig:sphereFlowAround}. 

The method is tested for several discretizations $h$ and regularizations $\epsilon$. We compare the computed drags to the theory by calculating the relative error of the drag in the $x$ direction as a percentage:

\[\text{Relative error (\%)} = \left (\frac{\lvert \tilde{f}_1 + 6 \pi \rvert}{6 \pi}\right ) \times 100 \] Results are tabulated in Table \ref{Tab:dragSphere}. Examining the table, we see that for a fixed $\epsilon$ ``small enough'' (in this case, for $\epsilon \leq 10^{-4}$), the relative errors decrease as the spatial discretization $h$ becomes smaller. For fixed $\epsilon \geq 10^{-2}$, the errors decrease as $h$ becomes smaller until $h \leq 0.1119$. Below this threshold, the errors appear to increase slightly and saturate.

\begin{table}[htbp!]
  \begin{subtable}[h]{0.9\textwidth}
    \centering
    \begin{tabular}{l|l|l?c|c|c|c|c|c}%
      \hline\hline & & & \multicolumn{6}{c}{\bfseries Regularization ($\epsilon$)}\\
      \bfseries no. tris & \bfseries DOF & \bfseries Disc. ($h$) & \bfseries $10^{-1}$ & $10^{-2}$ & $10^{-4}$ & $10^{-6}$ & $10^{-7}$ & $10^{-8}$
                                                                                                                                         \begin{filecontents*}{tables/relErrorDragSphere.csv}
numFaces,DOF,discretization,epsone,epstwo,epsthree,epsfour,epsfive,epssix
180,276,0.3675,1.9149,0.6850,0.9408,0.9434,0.9434,0.9434
720,1086,0.1860,2.5969,0.0103,0.2413,0.2438,0.2438,0.2438
2880,4326,0.0933,2.7792,0.1915,0.0597,0.0622,0.0622,0.0622
11520,17286,0.0467,2.8262,0.2380,0.0132,0.0157,0.0157,0.0157
\end{filecontents*}
\csvreader[head to column names]{tables/relErrorDragSphere.csv}{}
                                                                                                                                         {\\\hline  \numFaces & \DOF & \discretization & \epsone & \epstwo  & \epsthree & \epsfour & \epsfive & \epssix}
    \end{tabular}
    \caption{Relative drag force errors (\%) in x-direction for unit sphere translating with velocity $\vec U=(1,0,0)$.} 
    \label{Tab:dragSphere}
  \end{subtable} 

  \begin{subtable}[h]{0.9\textwidth}
    \centering
    \begin{tabular}{l|l|l?c|c|c|c|c|c}%
      \hline\hline & & & \multicolumn{6}{c}{\bfseries Regularization ($\epsilon$)}\\
      \bfseries no. tris & \bfseries DOF & \bfseries Disc. ($h$) & \bfseries $10^{-1}$ & $10^{-2}$ & $10^{-4}$ & $10^{-6}$ & $10^{-7}$ & $10^{-8}$
                                                                                                                                         \begin{filecontents*}{tables/relErrorRotationSphere.csv}
numFaces,DOF,discretization,epsone,epstwo,epsthree,epsfour,epsfive,epssix
180,276,0.3675,1.9149,0.6850,0.9408,0.9434,0.9434,0.9434
720,1086,0.1860,2.5969,0.0103,0.2413,0.2438,0.2438,0.2438
2880,4326,0.0933,2.7792,0.1915,0.0597,0.0622,0.0622,0.0622
11520,17286,0.0467,2.8262,0.2380,0.0132,0.0157,0.0157,0.0157
\end{filecontents*}
\csvreader[head to column names]{tables/relErrorRotationSphere.csv}{}
                                                                                                                                         {\\\hline  \numFaces & \DOF & \discretization & \epsone & \epstwo  & \epsthree & \epsfour & \epsfive & \epssix}
    \end{tabular}
    \caption{Relative net torque errors (\%) in z-direction for unit sphere prescribed to rotate with angular velocity $\vecg \Omega=(0,0,1)$.} 
    \label{Tab:momentSphere}
  \end{subtable}
  \caption{Net drag/moment errors for the a translating/rotating sphere.}
  \label{Tab:dragmom}
\end{table}

\begin{figure}[htbp!]
  \centering
  \begin{subfigure}{0.49\textwidth}
    \includegraphics[width=\textwidth]{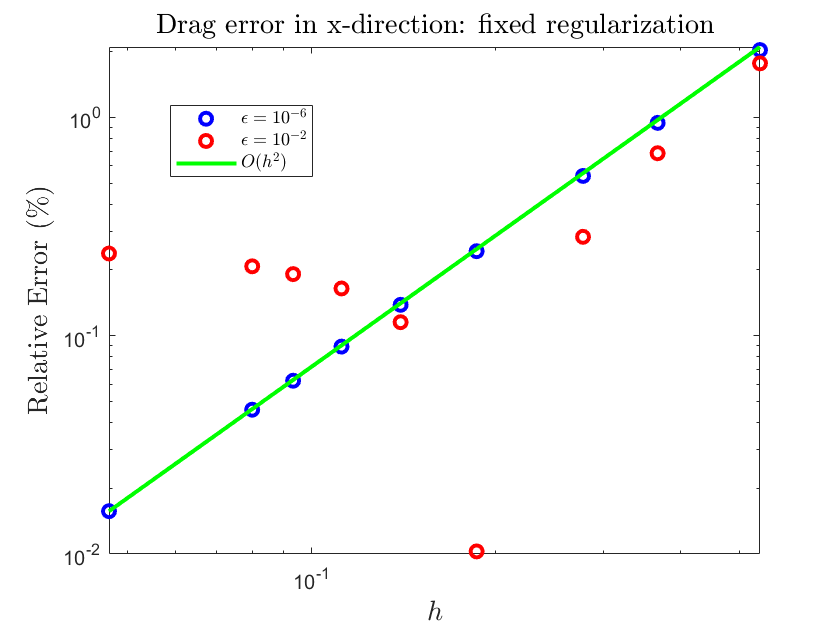}
    \caption{}
    \label{fig:dpa}
  \end{subfigure}
  \begin{subfigure}{0.49\textwidth}
    \includegraphics[width=\textwidth]{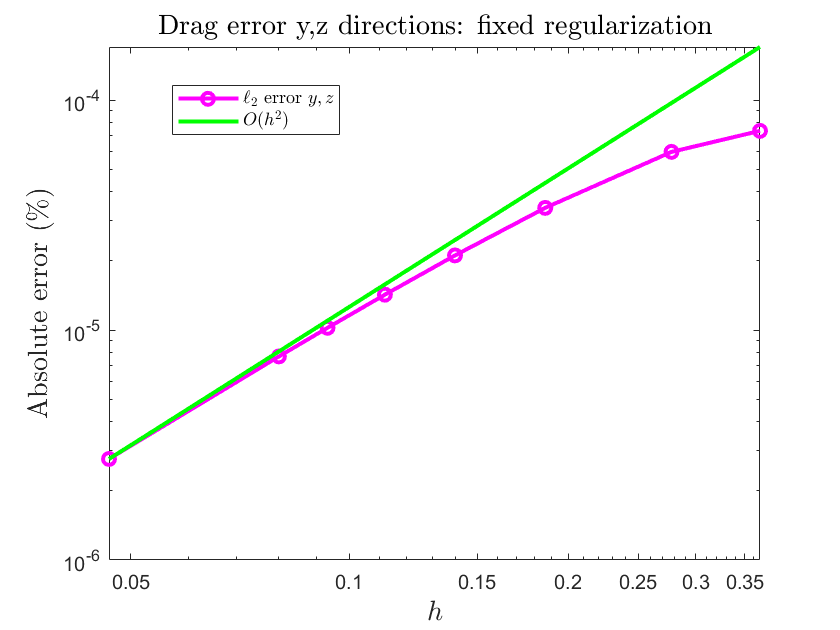}
    \caption{}
    \label{fig:dpb}
  \end{subfigure}
  \caption[Relative errors in the $x$-component of the drag for two regularizations]{Results from resistance problem for unit sphere prescribed to translate with velocity $\vec U=(1,0,0)$. (a) Relative errors in the $x$-component of the drag for two regularizations, $\epsilon=10^{-6}, 10^{-2}$. (b) Absolute errors in the $y,z$-components of the drag for $\epsilon=10^{-6}$.}
  \label{fig:dragProblem}
  \end{figure}

For $\epsilon \leq 10^{-4}$, the errors display uniform second order convergence in $h$. This is illustrated in Figure \ref{fig:dragProblem} for the regularization $\epsilon=10^{-6}$, with \ref{fig:dpa} showing the relative errors in the $x$-component of the drag, and \ref{fig:dpb} showing the absolute errors in the $y,z$ components of the drag. For comparison, we show the data for the regularization $\epsilon=10^{-2}$ in \ref{fig:dpa}.

\begin{figure}[h!]
  \centering
  \includegraphics[width=0.5\textwidth]{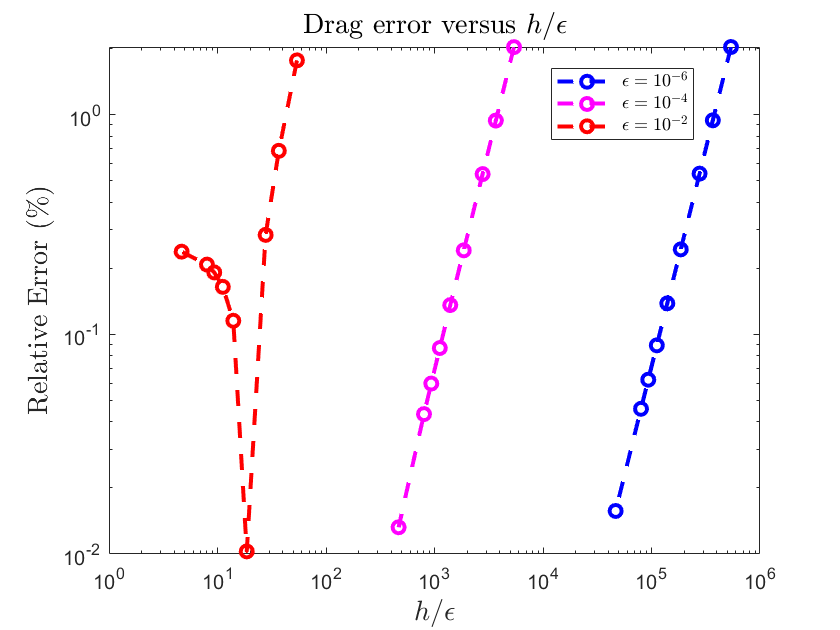}
  \caption[Drag errors as a function of $\epsilon/h$ for two fixed regularizations]{A comparison of the drag errors in the $x$-component as a function of $h/\epsilon$ for three regularizations: $\epsilon=10^{-2},10^{-4},10^{-6}$ . For the two smaller regularizations, the errors decrease for all $h$ tested whereas for the larger regularization ($\epsilon=10^{-2}$), the errors decrease until $h$ is about 10 times the size of $\epsilon$.}
  \label{fig:dragVsHovEps}
\end{figure}

From the data, it appears that as long as $\epsilon$ is chosen ``small enough,'' we will observe second order convergence in $h$ for a wide range of discretizations. As the red data in Figure \ref{fig:dragVsHovEps} show, there appears to be a limit to how small one can shrink $h$ relative to $\epsilon$ while maintaining gains in accuracy. These gains appear to reverse and then saturate once $h/\epsilon \approx 10$. For the magenta and blue curves, corresponding to $\epsilon = 10^{-4}$ and $\epsilon=10^{-6}$ respectively, the errors decrease at the same rate for all $h$ tested. \\[0.1em] 

\noindent \textit{Torque on a rotating sphere} \\
The rotational version of the resistance problem requires calculating the net torque, or moment, on a sphere of radius $a$ rotating uniformly with angular velocity $\vecg \Omega$. The theoretical value for the net torque is 

\[ \vec M = - 8 \pi \mu a^3 \vecg \Omega\]We again fix the parameters $a=1$ and the viscosity $\mu=1$. On the triangulated sphere, we prescribe the angular velocity $\vecg \Omega=(0,0,1)$ so that the sphere is rotating with unit speed about the $z$-axis. We calculate the approximate net moment $\widetilde{\vec M}=(\tilde{m}_1,\tilde{m}_2,\tilde{m}_3)$ by summing the torques from every triangle with equation \eqref{eq:netTorque} and comparing to the theoretical net torque $\vec M=(0,0,-8 \pi)$. 

\begin{figure}[htbp!]
  \centering
  \begin{subfigure}{0.49\textwidth}
    \includegraphics[width=\textwidth]{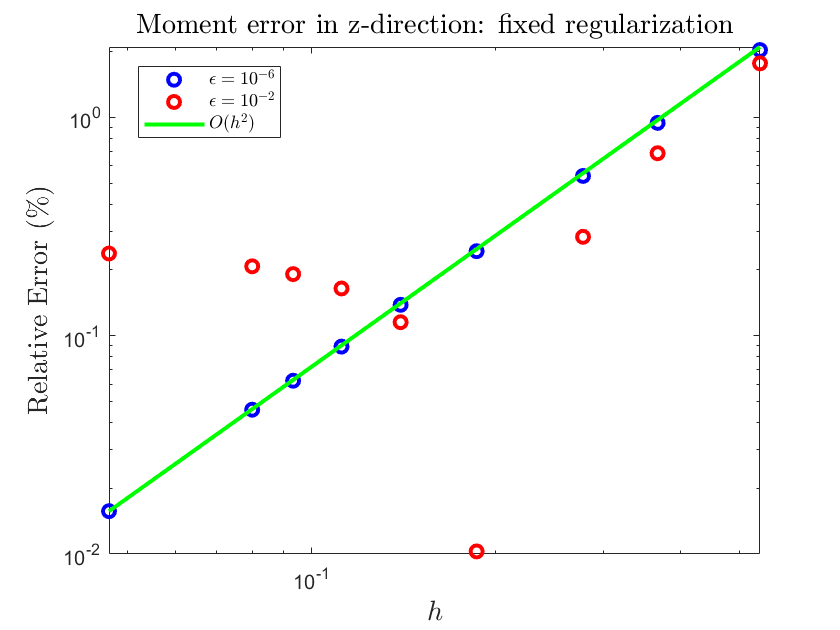}
    \caption{}
    \label{fig:mpa}
  \end{subfigure}
  \begin{subfigure}{0.49\textwidth}
    \includegraphics[width=\textwidth]{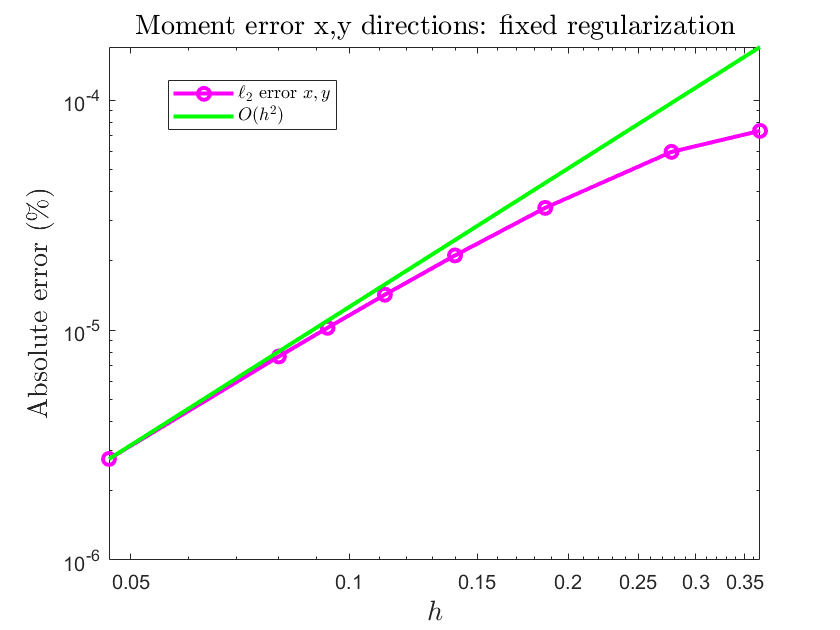}
    \caption{}
    \label{fig:mpb}
  \end{subfigure}
  \caption[Relative errors in the $z$-component of the net torque for two regularizations.]{Results from resistance problem for unit sphere prescribed to rotate with angular velocity $\vecg \Omega=(0,0,1)$. (a) Relative errors in the $z$-component of the net torque for two regularizations, $\epsilon=10^{-6}, 10^{-2}$. (b) Absolute errors in the $x,y$-components of the net torque for regularization $\epsilon=10^{-6}$.}
  \label{fig:momentProblem}
\end{figure}

We again calculate the relative error, which is now in the $z$-component of the net torque: 

\[\text{Relative error (\%)} = \left ( \frac{\lvert \tilde{m}_3 + 8 \pi \rvert}{8\pi}\right ) \times 100\]The results are tabulated in Table \ref{Tab:momentSphere}. The error behavior is similar to what was observed with the drag on the translating sphere. We see again that the errors are not uniformly second order in the spatial discretization $h$ until $\epsilon \leq 10^{-4}$. For the larger regularizations used, the error decays in $h$ for $h>0.1119$, but then begins to saturate for the smaller discretizations used. This contrast in error behavior is illustrated in Figure \ref{fig:mpa} for the regularizations $\epsilon=10^{-6}$ and $\epsilon=10^{-2}$. The absolute errors for the net torque in the $x,y$-components are shown in Figure \ref{fig:mpb} for regularization $\epsilon=10^{-6}$. Similar to what was seen in Figure \ref{fig:dpb}, the asymptotic decay begins when $h=0.2776$ and the errors are never larger than $10^{-4}$. \\[0.1em]

\noindent \textit{A note on condition numbers of resistance matrices}\\
We make one final observation about the resistance problem for the translating and rotating sphere. In the forward problem for the rotating sphere, we saw that the difference in performance between the regularized Stokeslet surfaces using a piecewise constant force density became comparable to the method when using a piecewise linear force density for the smallest $h$ used (compare last rows of Tables \ref{Tab:forwardLinearSurfacesTorque} and \ref{Tab:forwardConstantSurfacesTorque}.) However, when using the Stokeslet surfaces with piecewise constant force density in the resistance problem for the sphere, the condition numbers for the associated matrices are very large (Table \ref{Tab:condConstant}).

\begin{table}[ht]
\begin{subtable}[h]{0.8\textwidth}
 \begin{tabular}{l|l|l?c|c|c|c|c|c}%
    \hline\hline & & & \multicolumn{6}{c}{\bfseries Regularization ($\epsilon$)}\\
     \bfseries no. tris & \bfseries DOF & \bfseries Disc. ($h$) & \bfseries $10^{-1}$ & $10^{-2}$ & $10^{-4}$ & $10^{-6}$ & $10^{-7}$ & $10^{-8}$
                                                                                                                                                           \begin{filecontents*}{tables/condNumbersConstant.csv}
numFaces,DOF,discretization,epsone,epstwo,epsthree,epsfour,epsfive,epssix
180,540,0.3675,9.65e+16,2.58e+16,9.09e+14,1.11e+13,6.04e+11,1.77e+11
720,2160,0.1860,1.23e+17,1.29e+17,3.88e+16,7.81e+13,4.26e+12,2.85e+12
2880,8640,0.0933,7.34e+16,2.48e+17,4.22e+16,4.49e+14,6.20e+13,1.07e+13
11520,34560,0.0467,9.23e+17,1.21e+17,2.44e+17,4.27e+15,1.85e+14,2.01e+13
\end{filecontents*}
\csvreader[head to column names]{tables/condNumbersConstant.csv}{}
                                                                                                                                                  {\\\hline \numFaces & \DOF & \discretization & \epsone & \epstwo  & \epsthree & \epsfour & \epsfive & \epssix}
  \end{tabular}
  \caption{Condition numbers for matrices in resistance problem for translating/rotating sphere using piecewise constant forces.}
  \label{Tab:condConstant}
\end{subtable}
\begin{subtable}[h]{0.8\textwidth}
 \begin{tabular}{l|l|l?c|c|c|c|c|c}%
    \hline\hline & & & \multicolumn{6}{c}{\bfseries Regularization ($\epsilon$)}\\
     \bfseries no. tris & \bfseries DOF & \bfseries Disc. ($h$) & \bfseries $10^{-1}$ & $10^{-2}$ & $10^{-4}$ & $10^{-6}$ & $10^{-7}$ & $10^{-8}$
                                                                                                                                                           \begin{filecontents*}{tables/condNumbersLinear.csv}
numFaces,DOF,discretization,epsone,epstwo,epsthree,epsfour,epsfive,epssix
180,276,0.3675,1.36e+03,4.44e+02,3.90e+02,3.89e+02,3.89e+02,3.89e+02
720,1086,0.1860,3.83e+04,3.90e+03,3.02e+03,3.01e+03,3.01e+03,3.01e+03
2880,4326,0.0933,4.87e+06,3.98e+04,2.40e+04,2.39e+04,2.39e+04,2.39e+04
11520,17286,0.0467,1.41e+10,5.24e+05,1.93e+05,1.90e+05,1.90e+05,1.90e+05
\end{filecontents*}
\csvreader[head to column names]{tables/condNumbersLinear.csv}{}
                                                                                                                                                           {\\\hline \numFaces & \DOF & \discretization & \epsone & \epstwo  & \epsthree & \epsfour & \epsfive & \epssix}
  \end{tabular}
  \caption{Condition numbers for matrices in resistance problem for translating/rotating sphere using piecewise linear forces.}
\label{Tab:condLinear}
\end{subtable}
\label{Tab:cond}
\caption{Condition numbers for the resistance problem.}
\end{table}For the piecewise constant surfaces, we must take one point on each triangle (e.g. the centroid) and solve a $3N \times 3N$ system for the force density at each of the $N$ points, where $N$ is the number of triangles. This results in a force density that is discontinuous along the boundary of the triangles. For the piecewise linear surfaces, the size of the system depends on the number of unique points that are triangle vertices. In contrast to the piecewise constant surfaces, the resultant force density using piecewise linear surfaces is guaranteed to be continuous.

The condition numbers for the matrices in the translating/rotating sphere problem using piecewise linear surfaces (Table \ref{Tab:condLinear}) are several orders of magnitude smaller than the condition numbers from the matrices using piecewise constant surfaces. We conclude that there is a clear benefit in using piecewise linear surfaces over their piecewise constant counterparts in the inverse problem of computing the force density from a prescribed velocity.

\subsection{Rotating spheroid: Comparison of mesh design}
\begin{figure}[h!]  \centering
  \begin{subfigure}{0.42\textwidth}
{\includegraphics[width=\textwidth]{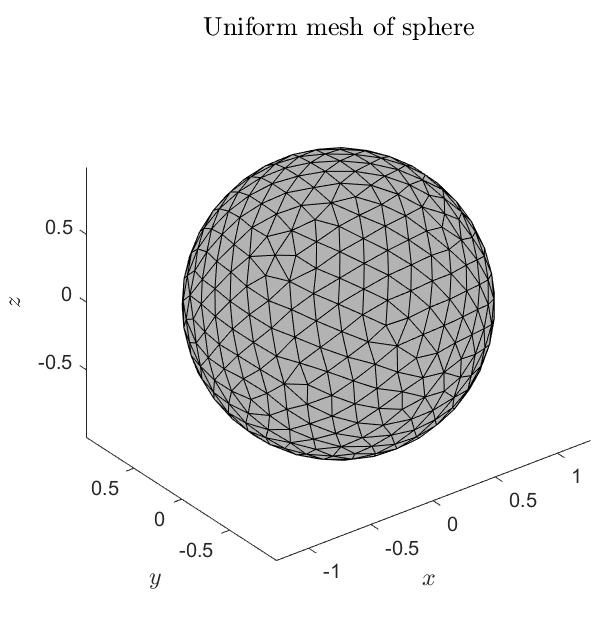}}
    \caption{}
    \label{fig:sphereuni}
  \end{subfigure}
  \begin{subfigure}{0.25\textwidth}
{\includegraphics[width=\textwidth]{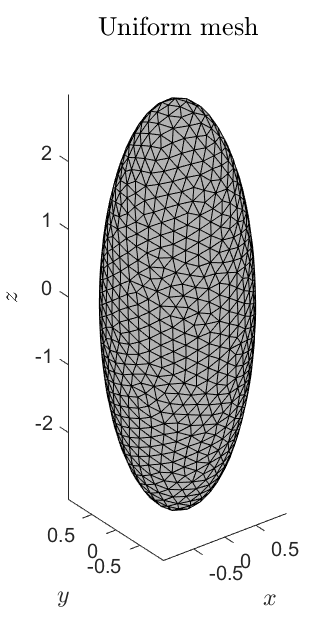}}
    \caption{}
    \label{fig:spheroiduni}
  \end{subfigure}
  \begin{subfigure}{0.27\textwidth}
{\includegraphics[width=\textwidth]{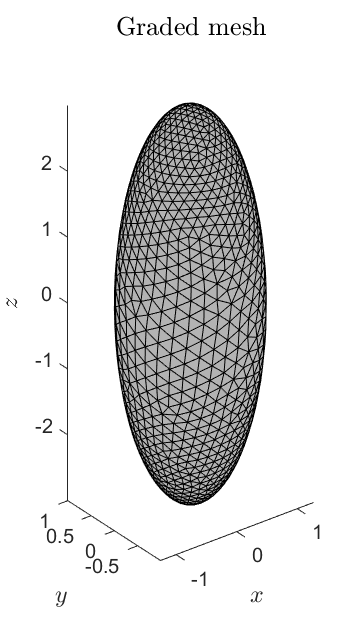}}
    \caption{}
    \label{fig:spheroidgraded}
  \end{subfigure}
  \caption{Meshes created using the Matlab package \textit{DistMesh}.
  (a) Uniform mesh of the unit sphere. (b/c) Uniform/graded mesh of the prolate spheroid (major:minor = 3:1) }
    \label{fig:mesh}
\end{figure}

    In this example, we turn to the problem of a prolate spheroid rotating about its major axis to investigate the implications of mesh design. We use the prolate spheroid described implicitly by $x^2+y^2+z^2/9=1$ so that the length of the major axis is $3$ and the length of the minor axis is $1$. To construct an approximately uniform mesh of the spheroid, we use the openly available Matlab package \textit{DistMesh} \cite{distmesh, persson2004}. We also use \textit{DistMesh} to construct a mesh of the unit sphere to compare results to the case of a sphere meshed with the same algorithm. The meshes mentioned in this discussion are shown in Figure \ref{fig:mesh}.

    
\begin{figure}[tb!]  \centering
  \begin{subfigure}{0.32\textwidth}
{\includegraphics[width=\textwidth]{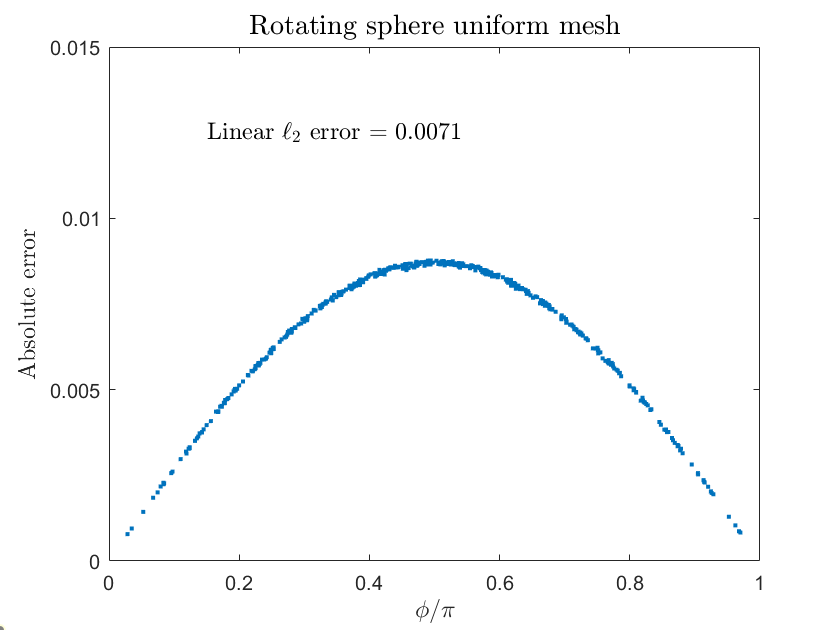}}
    \caption{}
    \label{fig:rotsphere}
  \end{subfigure}
  \begin{subfigure}{0.32\textwidth}
{\includegraphics[width=\textwidth]{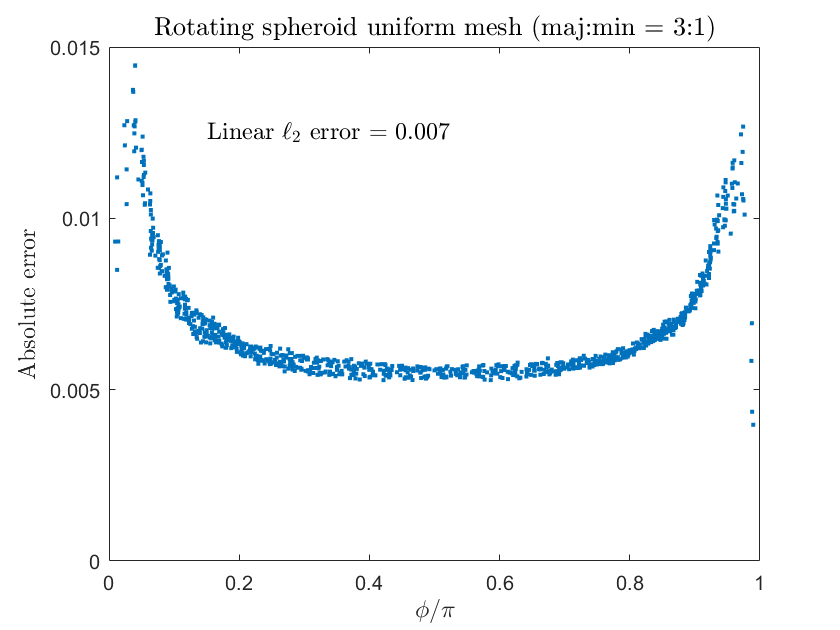}}
    \caption{}
    \label{fig:rotspheroiduni}
  \end{subfigure}
  \begin{subfigure}{0.32\textwidth}
{\includegraphics[width=\textwidth]{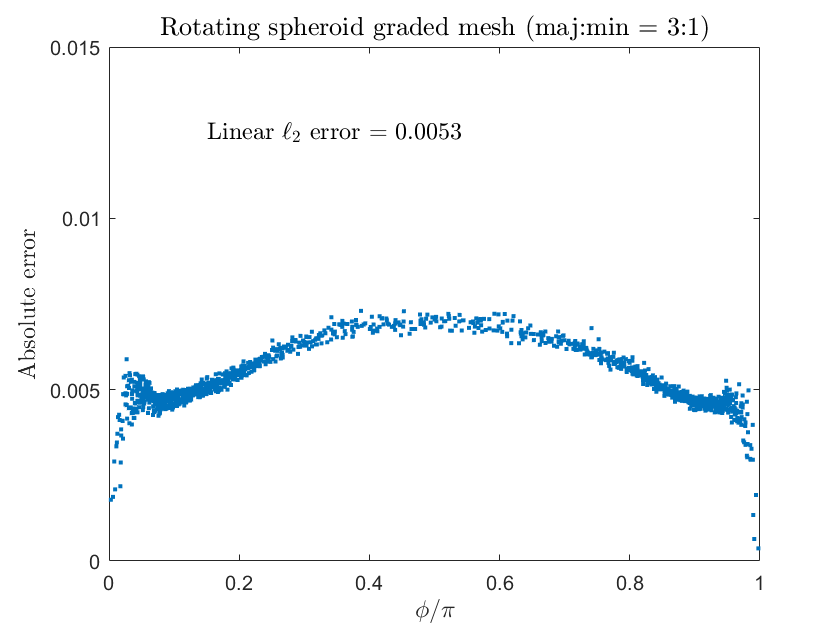}}
    \caption{}
    \label{fig:rotspheroidgraded}
  \end{subfigure}
  \caption{The effect of mesh design for the problem of a prolate spheroid rotating about its major axis at unit angular velocity. The results in (a),(b),(c) correspond to the meshes in (a),(b),(c) of Figure \ref{fig:mesh}. For a rotating sphere meshed uniformly (\ref{fig:rotsphere}), the errors are largest near the equator. For a prolate spheroid with major:minor axes in ratio 3:1 and meshed uniformly, the errors become largest near the poles (\ref{fig:rotspheroiduni}). If one uses a graded mesh that is finer near the poles, these errors decrease so that the shape of the error curve is similar to that of the sphere (\ref{fig:rotspheroidgraded}). There is still a slight increase near the poles, however, it is much smaller when compared to the uniform mesh.}
    \label{fig:rotspheroid}
\end{figure}

The hydrodynamic traction on a spheroid described by $z^2/a^2+r^2/b^2=1$ ($a\geq 1 > b$) rigidly rotating with unit angular speed about its major axis $\whvec e_z$ can be written \cite{kim2015},

    \begin{equation}
      \label{eq:tractionspheroid}
      \vec f(\vec x) = \frac{3(\whvec n(\vec x) \cdot \vec x)}{8\pi ab^4} \vec M \times \vec x 
    \end{equation}
    where $\vec x$ is a point on the spheroid, $\whvec n(\vec x)$ is the unit outward normal to the spheroid at $\vec x$, and $\vec M$ is the net torque exerted by the fluid on the spheroid. Using the notation of Chwang and Wu \cite{chwangwu74}, the net torque $\vec M$ is given by 

    \begin{equation*}
      \begin{split}
        \vec M &= - \frac{32}{3}\pi \mu a e \beta_0 \whvec e_z\\
        e \equiv \text{eccentricity}= \frac{\sqrt{a^2-b^2}}{a} &, \ \beta_0 = a^2 e^2 \left [ \frac{2e}{1-e^2} - \log \frac{1+e}{1-e} \right ]^{-1}         \end{split}
    \end{equation*}In the tests, we prescribe the traction at the spheroid vertices using \eqref{eq:tractionspheroid} with $(a,b)=(3,1)$, 

    Figure \ref{fig:rotspheroiduni} displays the pointwise errors at the vertices as a function of the scaled polar angle $\phi/\pi$. Comparing to Figure \ref{fig:rotsphere} which shows the errors for the case of a sphere, we notice a difference in the shape of the error curves. In particular, errors near the equator of the ellipsoid are smaller relative to the errors near the poles. This contrasts to the errors on the sphere which were smallest near the poles.

    In order to reduce the errors near the poles, one can use a graded mesh which uses smaller triangles in their vicinity (see Figure \ref{fig:spheroidgraded}.) Results for this case are plotted in Figure \ref{fig:rotspheroidgraded}. We notice an improvement in accuracy near the poles resulting in an error curve that more closely resembles that seen for the sphere. There is still a small region where the errors appear to slightly increase near the poles. However, these errors eventually do decrease even closer to the poles. This example illustrates the dependence of the error on the particular mesh design. For objects like the spheroid that possess axes of symmetry of different lengths, it may be advantageous to use nonuniform meshes. 

\subsection{The squirmer model}
The ``squirmer'' is a canonical model of swimming at zero Reynolds number, first introduced by Lighthill \cite{lighthill52} and later extended to a deforming envelope model by Blake \cite{blake71}. We focus on the former, which is a simpler model of ciliate swimming when compared to the latter.

The swimming organism is modeled as a sphere which propels itself with velocity $\vec U$ by generating a tangential slip velocity $\vec u_{\text{slip}}(\vec x)$ at points $\vec x$ on the sphere. From the modeling perspective, the slip velocity can be understood as the movement of the beating cilia tips. This coordinated beating generates the viscous propulsive force which balances the hydrodynamic drag experienced by the squirmer.

Mathematically, the problem is to solve the Stokes equations with the condition that points on the surface of the body, $\partial B$, undergo a rigid translation and rotation with velocity $\vec U$ and angular velocity $\vecg \Omega$ respectively, in addition to the imposed slip velocity $\vec u_{\text{slip}}(\vec x)$. This boundary condition is stated mathematically as  

\[\vec u(\vec x) \large \rvert_{\partial B} = \vec U + \vecg \Omega \times (\vec x - \vec X_{\text{center}})+ \vec u_{\text{slip}}(\vec x)  \]where $\vec X_{\text{center}}$ is the sphere center. Here, $\vec U$ and $\vecg \Omega$ are unknowns which must be solved for. This creates 6 additional unknowns that necessitate two more equations to close the system. The physically relevant choice is to enforce the net force and torque on the organism to be zero. In other words, the organism experiences no external force or torque so that its movement is entirely due to its coupling with the fluid.

\begin{figure}[htbp!]
    \centering
    \begin{subfigure}{0.49\textwidth}
      \includegraphics[width=\textwidth]{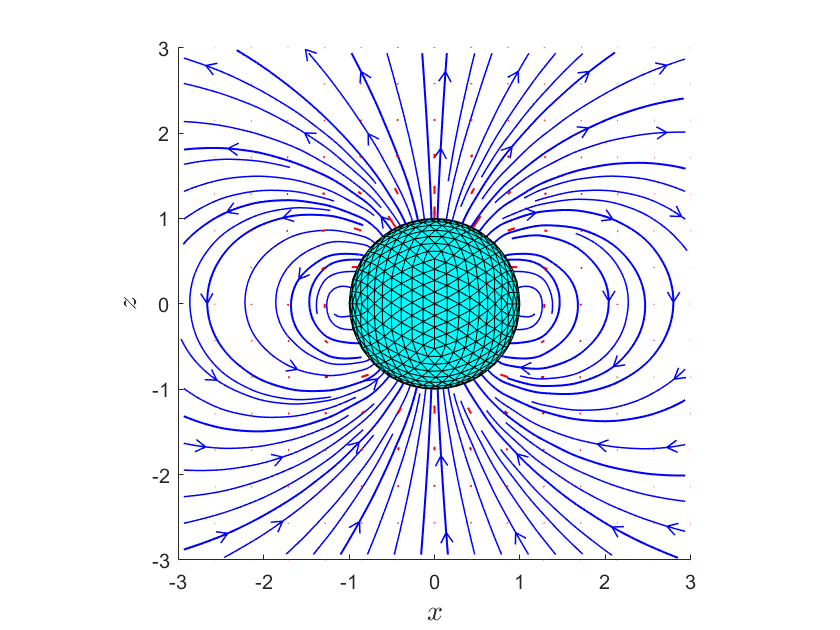}
      \caption{Streamlines due to squirmer in lab frame}
    \end{subfigure}
    \begin{subfigure}{0.49\textwidth}
      \includegraphics[width=\textwidth]{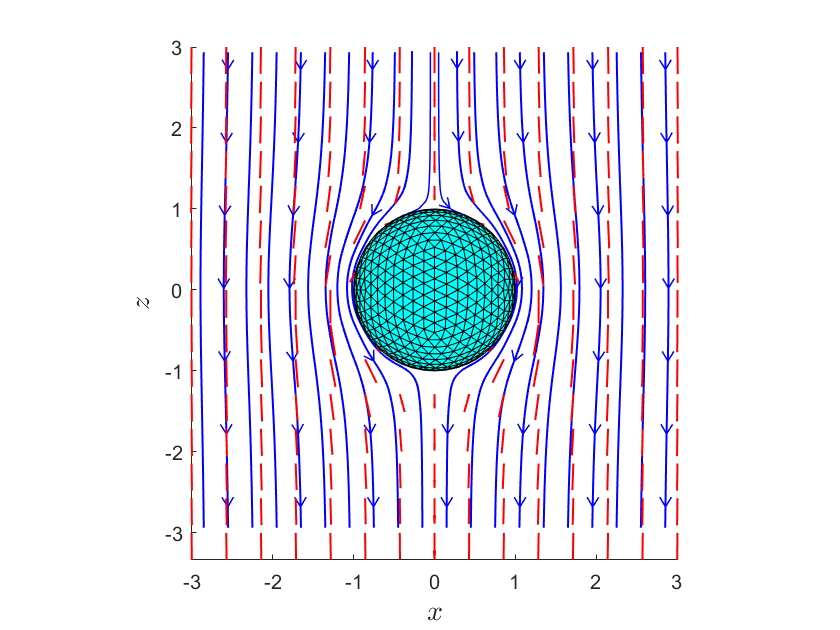}
      \caption{Streamlines due to squirmer in co-moving frame}
    \end{subfigure}
  \caption{Flow visualization near squirmer of radius $a=1$ using average spatial discretization $h=0.1398$ and regularization $\epsilon=10^{-4}$.}
  \label{fig:squirmerStreamlines}
\end{figure}

The complete statement of the problem is to solve for the velocity field $\vec u$, the translational velocity $\vec U$, and the angular velocity $\vecg \Omega$ such that 

\begin{equation}
  \label{eq:swimproblem}
  \begin{split}
    \vec 0 &= -\nabla p + \mu \Delta \vec u  \\
    0 &= \nabla \cdot \vec u \\
    \vec 0 &= \int \int_{\partial B} \vec f(\vec s) \ d\vec s  \\
    \vec 0 &= \int \int_{\partial B} (\vec x(\vec s) - \vec X_{\text{center}}) \times \vec f(\vec s) \ d \vec s \\
    \vec u(\vec x) \large \rvert_{\partial B} &= \vec U + \vecg \Omega \times (\vec x - \vec X_{\text{center}})+ \vec u_{\text{slip}}(\vec x) 
  \end{split}
\end{equation}The solution, first derived by Lighthill, proceeds by imposing an axisymmetric slip velocity of the form

\[\vec u_{\text{slip}}(\theta) = \begin{bmatrix} u_r(\theta) \\ u_\theta(\theta) \\ u_\phi(\theta) \end{bmatrix}= \begin{bmatrix} 0 \\ \sum_{n=1}^{\infty}B_n V_n(\cos \theta) \\ 0 \end{bmatrix}\]where $\theta \in [0,\pi]$ is the angle from the north pole of the sphere, $u_r,u_\theta,u_\phi$ are the radial, polar, and azimuthal components of the velocity field respectively, $B_n$ are coefficients, and $V_n(x)$ is a function defined in terms of the Legendre polynomial of degree $n$, denoted $P_n(x)$: 

\[V_n(x)=\frac{2 \sqrt{1-x^2}}{n(n+1)}P_n'(x) \]The resultant swimming velocity, $\vec U$, is determined by the first coefficient $B_1$ and is in the direction of the north pole $\vec n$ \cite{ishimoto21}:

\[\vec U = \frac 2 3 B_1 \vec n\]The higher order modes determine the decay of the flow far away from the squirmer but do not contribute to the swimming velocity. In the following, we set $B_1=3/2$ and $B_n = 0 $ for $n \geq 2$.With these parameters, the squirmer swims at unit speed in the direction of the north pole $\vec n$, and its angular velocity is $\vecg \Omega = \vec 0$. In spherical coordinates with origin at the center of the sphere, the solution for the velocity field  $\vec u$ is given in terms of the radial component of the velocity $u_r(r,\theta)$ and the polar component of the velocity $u_\theta(r,\theta)$. The azimuthal component of the velocity field, $u_\phi(r,\theta)$, is uniformly zero:

\begin{equation}
  \label{eq:squirmer}
  \begin{split}
    u_r(r,\theta) &= \left (\frac{a}{r}\right)^3 \cos \theta \\
    u_\theta(r,\theta) &= \frac{1}{2} \left (\frac{a}{r} \right )^3 \sin \theta \\
    u_\phi &= 0
  \end{split}
\end{equation}To use the method of Stokeslet surfaces, we first triangulate
the sphere using the icosahedral mapping described in the previous section. Using the imposed slip velocity,

\[\vec u_{\text{slip}}(\theta) = \begin{bmatrix} 0 \\ \frac{3}{2}V_1(\cos \theta) \\ 0 \end{bmatrix}\]we solve the discrete version of the system of equations in \eqref{eq:swimproblem} for the force density $\vec f$, the translational swimming velocity $\vec U$, and the angular velocity $\vecg \Omega$. Figure \ref{fig:squirmerStreamlines} shows streamlines for the numerical solution of this flow in the lab frame and in a frame moving with the organism. 

  
\begin{figure}[ht!]
  \centering
  \begin{subfigure}{0.49\textwidth}
    \includegraphics[width=\textwidth]{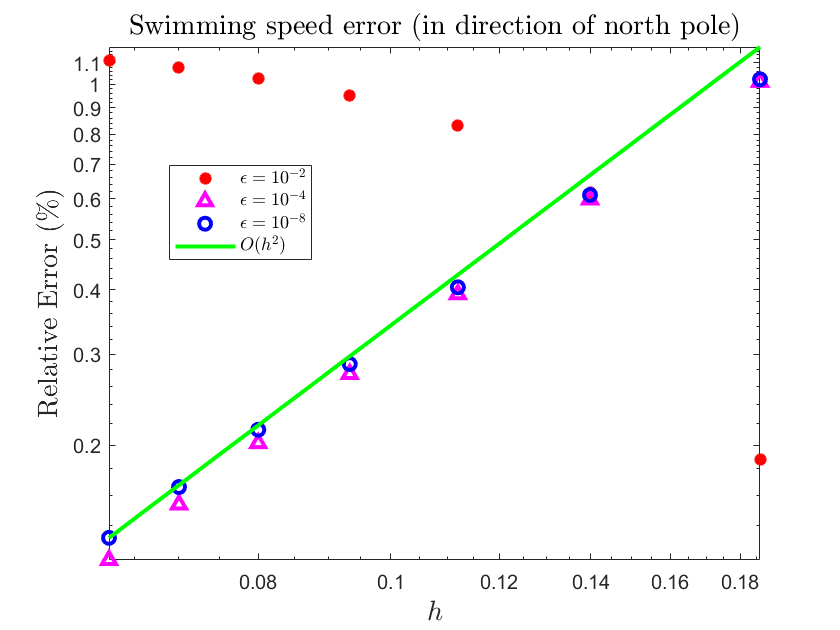}
    \caption{}
    \label{fig:esqa}
  \end{subfigure}
  \begin{subfigure}{0.49\textwidth}
    \includegraphics[width=\textwidth]{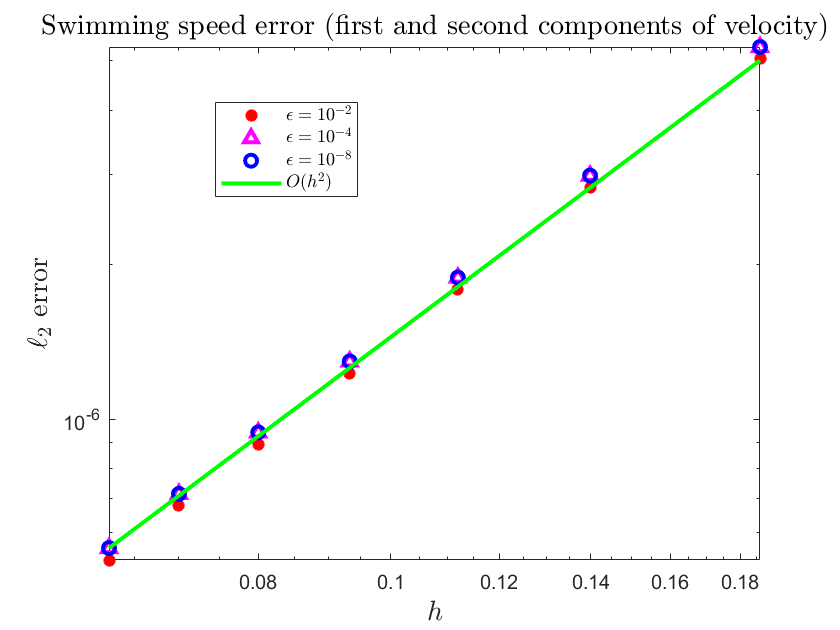}
    \caption{}
    \label{fig:esqb}
  \end{subfigure}\\
  \begin{subfigure}{0.49\textwidth}
    \includegraphics[width=\textwidth]{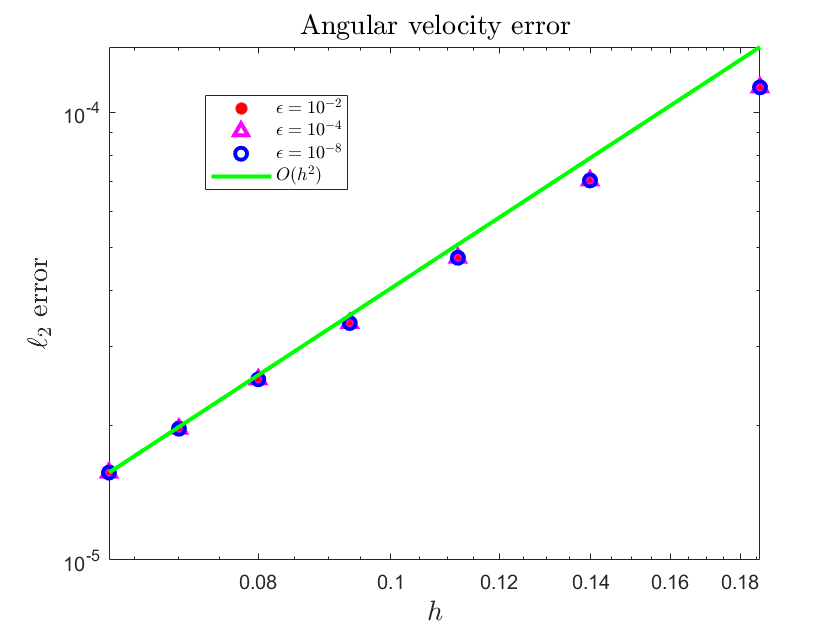}
    \caption{}
    \label{fig:esqc}
  \end{subfigure}
  \caption{(a) Relative errors in third component of swimming velocity $\vec U$. (b) Absolute errors in the first and second components of the swimming velocity $\vec U$. (c) Absolute error in the squirmer angular swimming velocity $\vecg \Omega$.}
  \label{fig:errorSquirmer}
\end{figure}

We tested the model using several spatial discretizations of the squirmer ($h$) and regularizations $\epsilon=10^{-2},10^{-4},10^{-6}$. We compare the numerical results to the theory by computing the relative error in the third component of the translational velocity $\vec U$ and the absolute error for the first and second components. As shown in Figure \ref{fig:esqa}, the relative error in the third component decays at a second order rate in $h$ as long as the regularization is chosen sufficiently small relative to $h$. If $\epsilon$ is too large relative to $h$, then this convergence is not observed: the errors increase but appear to saturate as $h$ gets smaller (red data in \ref{fig:esqa}.) We note that this is the same error behavior observed in the previous example of the translating/rotating sphere for the larger $\epsilon$ tested. 

We also compute the $\ell_2$ errors in the first and second components of the swimming velocity $\vec U$. Since the theoretical swimming velocity is only in the $z$ direction, the error computed is the magnitude of the swimming speed in all other directions besides $\whvec n$. As expected, these errors are small ($\sim 10^{-6}$) and decay at a second order rate in $h$ for all regularizations (Figure \ref{fig:esqb}.) In Figure \ref{fig:esqc}, we plot the errors in the squirmer's angular velocity, $\vecg \Omega$. Since the theoretical angular velocity of the squirmer is $\vecg \Omega=\vec 0$, we plot the absolute errors which are typically on the order of $10^{-5}$ indicating small deviations from the theoretical angular swimming velocity. For these measurements, second order convergence is observed for all regularizations tested.

To illustrate the gains in accuracy achieved with the analytic integration as opposed to a method that relies on quadrature, we compare the near-field fluid velocity error from Stokeslet surfaces with an implementation of the boundary element method of regularized Stokeslets (BEMRS\cite{smith2009}). In our implementation, the boundary elements are the triangles and constant forces are used over each element. The velocity is prescribed at the centroid of each element so that the number of unknowns in the system is $3(N+2)$ where $N$ is the number of triangles, and the additional $2$ comes from the zero net force/torque condition. Like the examples in \cite{smith2009}, we use a relatively expensive Gaussian quadrature for the nearly singular integrals (a nonproduct rule of degree 23 requiring 102 points in the triangle \cite{quad102}) and a less expensive Gaussian quadrature (a nononproduct rule of degree 7 requiring 15 points in the triangle \cite{quad15}) for the other integrals. We consider an integral to be nearly singular when the centroid and the Gaussian quadrature points lie on the same triangle. We used the openly available Matlab package \textit{quadtriangle} written by Ethan Kubatko \cite{quadtriangle} to implement the quadrature rules.

The sphere surface is discretized by 1280 triangles, resulting in an average spatial discretization $h \approx 0.1398$. For the Stokeslet surfaces, this corresponds to solving for the force density at $N_1=642$ points, while for the boundary element regularized Stokeslet method, we solve for the force density at each centroid, so we have $N_2=1280$ points that correspond to unkowns in the system. For each method, we solve a linear system ($1932 \times 1932$ system for the Stokeslet surfaces, $3846 \times 3846$ system for the BEMRS) for the force density on the squirmer surface. The regularization is fixed at $\epsilon=10^{-4}$ for both methods. Once the force density is solved for, the velocity field can be evaluated at any other point using \eqref{eq:velotri} for the Stokeslet surfaces, and a quadrature method for the BEMRS. For the BEMRS evaluation, we use the 15 point Gaussian quadrature method \cite{quad15} mentioned earlier.

\begin{figure}[tb!]
  \centering
  \begin{subfigure}{0.49\textwidth}
    \includegraphics[width=\textwidth]{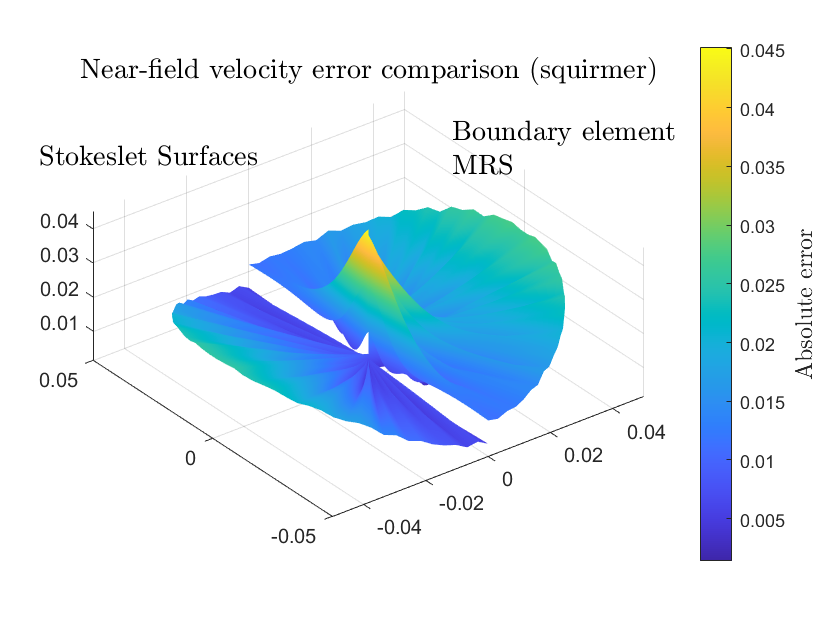}
    \caption{}
    \label{fig:surfvsbe1}
  \end{subfigure}
  \begin{subfigure}{0.49\textwidth}
    \includegraphics[width=\textwidth]{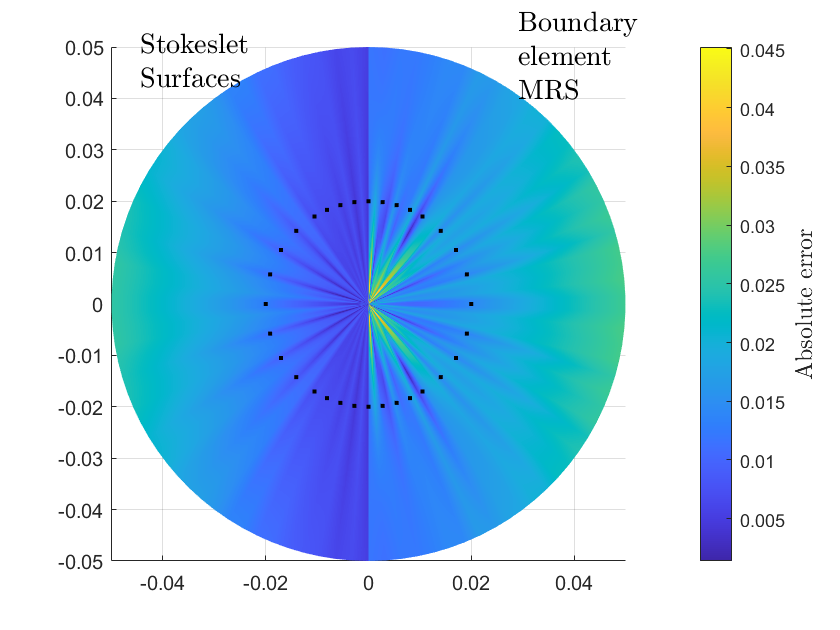}
    \caption{}
    \label{fig:surfvsbe2}    
  \end{subfigure}
  \begin{subfigure}{0.49\textwidth}
    \hspace{3em}
    \scalebox{0.9}{
    \def\svgwidth{\columnwidth}
    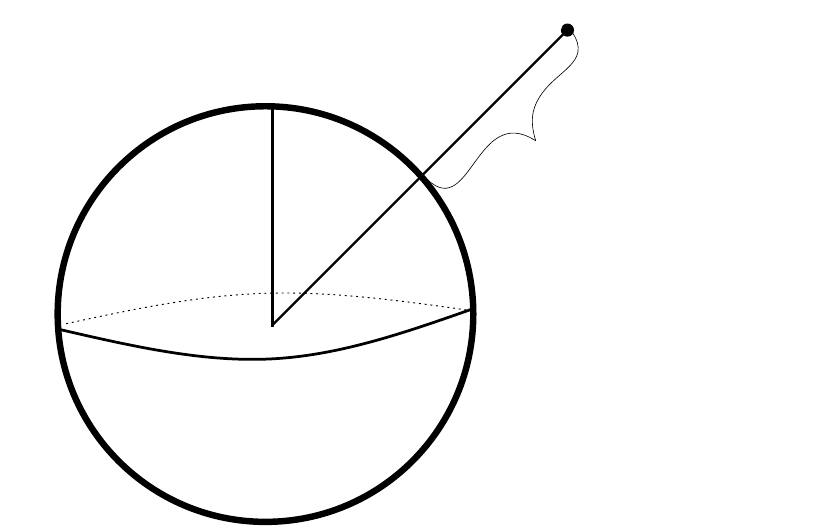
}
    \caption{}
    \label{fig:surfvsbe3}
\end{subfigure}
  \caption[Surfaces error comparsion to boundary element method]{(a) Absolute error in the $y=0$ plane as a function of $r$ and $\phi$, where $r$ is the distance from the unit sphere (squirmer). (b) A flattened version of the surface plot. The black dots drawn correspond to triangle vertices in the $y=0$ plane with the same polar angle. The smallest errors occur in line with these vertices for the Stokeslet surfaces. (c) In (a) and (b), the data points are plotted in polar coordinates as $r$ and $\phi$. $r$ corresponds to the distance from the squirmer surface and $\phi$ corresponds to the regular polar angle as shown in the figure.}
  \label{fig:surfvsbe}
\end{figure}

In Figure \ref{fig:surfvsbe1} and \ref{fig:surfvsbe2}, the near-field fluid velocity error is plotted as a surface and color coded by the pointwise absolute error. The radial distance from the origin corresponds to the distance from the evaluation point to the surface of the unit sphere (data points closer to the origin are physically closer to the sphere) and the polar angle corresponds to the regular polar angle $\phi$ in spherical coordinates (see Figure \ref{fig:surfvsbe3}.) The errors shown are evaluated in the $y=0$ plane of physical space and we evaluate the velocity at points up to $0.05$ units from the surface of the unit sphere (5 \% of the squirmer radius.) Given the symmetry of the problem and the  discretization of the sphere, we show the errors for the Stokeslet surfaces on the left and the errors for the BEMRS implementation on the right.

Examining Figure \ref{fig:surfvsbe1}, we observe that as the evaluation point gets closer to the sphere surfaces, the pointwise $\ell_2$ error for the BEMRS peaks at about $0.045$. However, the errors closest to the surface of the sphere are minimized for the Stokeslet surfaces. Moreover, the darker shades of blue indicate that the errors at corresponding points are smaller for Stokeslet surfaces compared to those of the BEMRS.

As the evaluation point moves further away from the sphere, the errors become comparable. Both methods have similar error magnitudes near the polar angle $\phi = \pm \pi/2$ farther away from the sphere surface (about $r=0.04$ from the surface), however, the region where these errors occur is larger for the BEMRS. In Figure \ref{fig:surfvsbe2}, one can see ``folds'' along different polar angles $\phi$ where the error appears to oscillate. These valleys in the error correspond to triangle vertices in the $y=0$ plane and these vertices are represented by black dots corresponding to their polar angle $\phi$.

This example illustrates the difficulty in using traditional quadrature methods for nearly singular integrals. On the other hand, it demonstrates the power of the analytic integration inherent in our  method.

\subsection{Internal flow in a pipe}
As a final example, we simulate internal flow in a pipe with rectangular cross section. We take the centerline of the pipe to be along the x-axis and $-a<y<a$, $-b<z<b$ so that the dimensions of a cross section are $2a \times 2b$. The fluid is driven by a constant pressure gradient in the x-direction ($\pp{p}{x}=-\Delta P$) and we set the boundary velocity $\vec u \rvert_{y=\pm a} = \vec u \rvert_{z=\pm b}=\vec 0$. The analytical solution is $\vec u_{\text{pipe}} = (u(y,z),0,0)$ where (\cite{pozr2011})

\begin{equation}
  \label{eq:uPoiseuille}
  \begin{split}
    u(y,z)&=\frac{\Delta P}{2 \mu}\left [b^2 -z^2 + 4b^2 \sum_{n=1}^{\infty}\frac{(-1)^n}{\alpha_n^3}\frac{\cosh(\alpha_n \frac{y}{b})}{\cosh(\alpha_n \frac{a}{b})}\cos(\alpha_n \frac{z}{b}) \right ] \\
          &\alpha_n \equiv (n-\frac 1 2)\pi \ 
  \end{split}
\end{equation}

\begin{figure}[h!]
  \centering
  \scalebox{0.6}{
    \def\svgwidth{\columnwidth}
    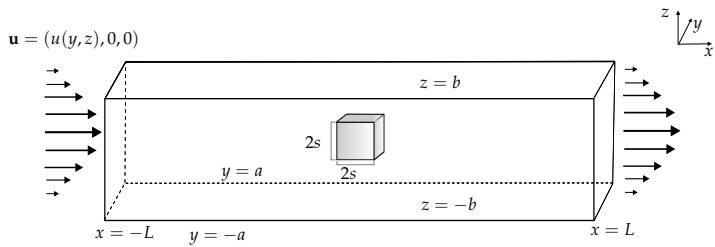
}
  \caption{Problem setup for the obstructed flow in a rectangular pipe}
  \label{fig:internalFlowGeometry}
\end{figure}

We will focus on the specific case of a square cross section so that $a=b$. To test the Stokeslet surfaces, we place a cube of side length $2s$ in the pipe as shown in Figure \ref{fig:internalFlowGeometry}. One end of the pipe is at $x=-L$ and the other end is at $x=L$. The cube is placed so that its center lies at the origin of the domain.  

The faces of the pipe and the cube are triangulated. To simulate the background flow $\vec u_{\text{pipe}}$, we first solve a linear system for the forces such that the cube velocity $\vec U_{\text{cube}}=-\vec u_{\text{pipe}}$ at the cube nodes and the pipe velocity $\vec U_{\text{pipe}}=\vec 0$ at the pipe nodes. Then we add the background flow $\vec u_{\text{pipe}}$ so that the velocity at the cube nodes is canceled out and the cube is stationary.

We choose the parameters for the simulation to be $a=b=1$, $\Delta P=1$, $\mu=1$, $s=0.25,$ and $L=2.5$. Since $2L/2s=10$, the length of the pipe is 10 times the side length of the cube and the effects of the obstruction are small at the inlet and outlet of the pipe so the flow is approximately the prescribed background velocity at $x=\pm L$. Since the analytical solution \eqref{eq:uPoiseuille} is an infinite series, we truncate the sum at $N=50$ terms (as the series is alternating, the error from the truncation can be estimated as being no larger than $\frac{\Delta P}{2\mu}4b^2$ times the maximum value of the next term in the series: $\approx 5 \times 10^{-7}$ using our problem parameters).

To test how well the method prevents fluid from leaking through solid boundaries, we measure the absolute flux through the front face of the cube. Here, we define the absolute flux as 

\begin{equation}
  \label{eq:absflux}
  \text{absolute flux} = \lvert \text{flux} \rvert = \int_F \lvert \vec u \cdot \whvec n \rvert \ d\vec s
\end{equation}where $F$ is the face of an oriented surface in the fluid. Using the absolute flux means we do not distinguish between fluid going in or out of a body. We approximate the absolute flux by using a $4\times 4$ Gaussian quadrature product rule over each triangle and summing the contributions.

We measure the leak by comparing the absolute flux over the front face of the cube with the absolute flux over the same surface if the pipe were not obstructed. For this particular problem, the flux through the front face of the cube would be 

\begin{equation}
  \label{eq:fluxNoCube}
  \begin{split}
    \lvert \text{flux without cube}\rvert &=  \int_{-s}^s \int_{-s}^s u(y,z) \ dy \ dz \\
                                          &=-\frac 2 3 \frac{\Delta P}{\mu} s^2 \left ( -3b^2 + s^2\right ) \\
    &+ \frac{\Delta P}{2\mu} 4b^4 \sum_{n=1}^{\infty} \frac{(-1)^n}{\alpha_n^5}\frac{\sinh\left(\alpha_n \frac{s}{b}\right)}{\cosh\left(\alpha_n \frac{a}{b}\right)}\sin\left(\alpha_n \frac{s}{b}\right)
  \end{split}  
\end{equation}The leak is then defined,

\begin{equation}
  \label{eq:leak}
  \text{leak} = \frac{\lvert \text{flux}\rvert}{\lvert \text{flux without cube}\rvert}
\end{equation}Like in the previous examples, we test the method with different combinations of spatial discretizations and regularizations. The surface of the cube is triangulated by placing a grid of length $h_{\text{cube}}$ over each face and drawing a diagonal from one corner to another, creating two right triangles. The surface of the rectangular pipe is triangulated in the same way but with another grid size $h_{\text{pipe}}$. For all the tests, we fix $h_{\text{pipe}}=0.2$ and vary $h_{\text{cube}}$ from $0.02$ to $0.1$. Since $h_{\text{pipe}}$ is fixed for all tests, we will refer to $h_{\text{cube}}$ as $h$ in the following analysis.

\begin{figure}[htbp!]
  \centering
  \begin{subfigure}{0.49\textwidth}
    \includegraphics[width=\textwidth]{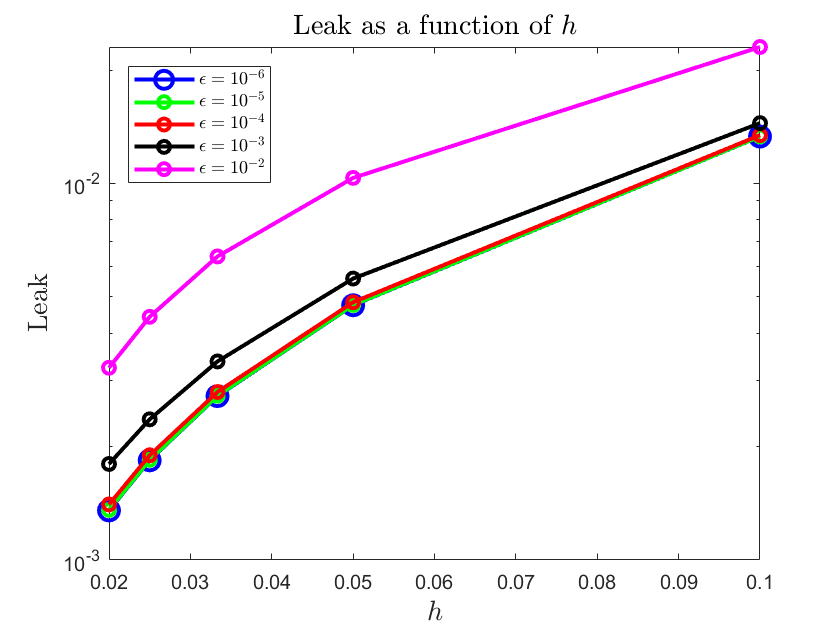}
    \caption{}
    \label{fig:leaka}
  \end{subfigure}
  \begin{subfigure}{0.49\textwidth}
    \includegraphics[width=\textwidth]{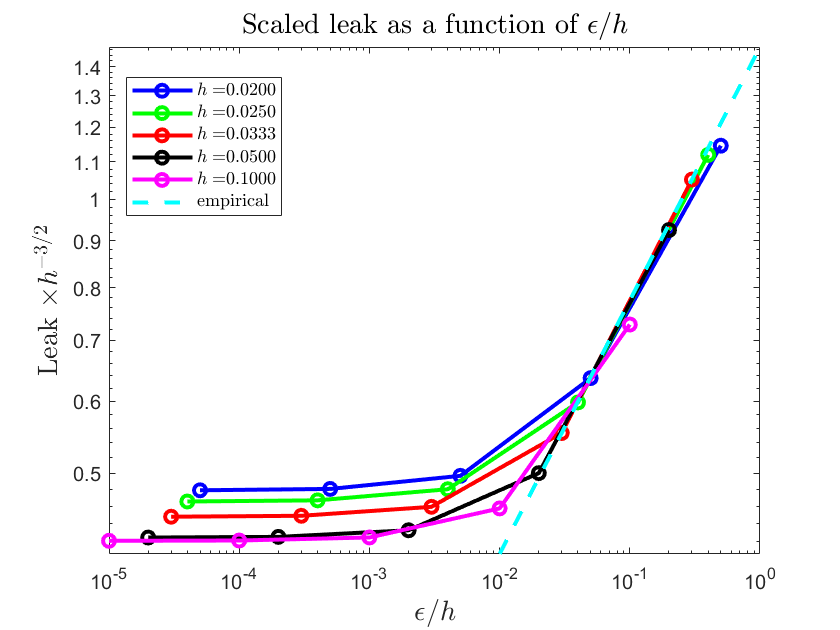}
    \caption{}
    \label{fig:leakb}    
  \end{subfigure}
  \caption[Leak as a function of spatial discretization $h$ in internal flow]{As $\epsilon$ and $h$ are decreased, the leak through the front face is reduced. However, the reduction in error that comes from reducing $\epsilon$ appears to saturate when $\epsilon \approx 10^{-4}$.}
  \label{fig:leak}
\end{figure}

In Figure \ref{fig:leaka}, the leak is shown as a function of the cube discretization $h$ with different regularizations. For all fixed regularizations, the leak is reduced as $h$ decreases. On the other hand, for a fixed discretization $h$, the returns in leak reduction from making $\epsilon$ smaller diminish once $\epsilon \approx 10^{-4}$.

In Figure \ref{fig:leakb}, the leak curves are scaled by $h^{-3/2}$ and plotted as a function of $\epsilon/h$. Now each curve corresponds to a fixed spatial discretization of the cube. Here, we see that the scaled leaks fall approximately on the same curve. When $10^{-2}<\epsilon/h<1$, we observe exponential decay in the scaled leaks.  When $\epsilon/h<10^{-2}$, all the curves flatten out. Due to the logarithmic scaling, the difference between the curves appears largest in this region. In reality, the difference between the any two curves is less than $10^{-1}$. The cyan colored curve was created from the results for the median discretization, $h=0.0333$, using the data for $10^{-2}<\epsilon/h<1$. This gives us the approximation

\[\frac{\text{Leak}(\epsilon/h)}{h^{3/2}} = 1.4709 \left(\frac{\epsilon}{h}\right)^{-0.2788} \]This empirical result is similar to the one given in \cite{cortez18} for the leak using regularized Stokeslet segments. In that case, however, the leak was defined in terms of the total error, not just the flux through the immersed body like was done here. Also, the leak in \cite{cortez18} was scaled differently (by $h^{1/2}$) and was shown to decrease as $\epsilon/h$ increased. In contrast, the scaled leak for the surfaces decreases as a function of $\epsilon/h$.

\section{Conclusions}
The variation on the method of regularized Stokeslets (MRS) introduced in this work was motivated by weakening the dependence of the choice of regularization parameter, $\epsilon$, on the spatial discretization, $h$, in the classic implementation of MRS. The usual choice of $\epsilon = O(h)$ is justified by the quadrature error estimates shown in \cite{cortez2005} which take the asymptotic form $O(h^p/\epsilon^q)$ with $p<q$. Using $\epsilon$ much smaller than $h$ yields large errors in the space between discretization nodes, resulting in fluid leaking through surfaces.  

The approach introduced in this work removes the need for quadrature by replacing it with exact integration over a triangulated surface. In effect, one can choose $\epsilon$ to be several orders of magnitude smaller than $h$ and observe second order convergence in $h$ without having to tune $\epsilon$. 

As a matter of context, it is important to recognize that an advantage of the typical regularized Stokeslet methods compared to other methods (e.g. the boundary element method) is its ease of implementation. This is due in part to the meshless nature of most implementations. This work is a departure from that paradigm, which may be better suited for some applications. On the one hand, the accurate modeling of slender bodies like a flagellum  may not require a surface meshing and is more efficiently implemented using the classical MRS or regularized Stokeslet segments \cite{cortez18}. On the other hand, in the case of large boundaries (e.g. a tube or channel), it may be more costly to discretize the surface by points. A surface meshing, in conjunction with the regularized Stokeslet surfaces method like the one presented here, allows the practitioner to use far fewer points, resulting in fewer interaction term evaluations and a smaller linear system to solve.

In several examples, the expected second order convergence in $h$ was numerically validated. These examples were also illustrative of instances where some care must be taken. For instance, we saw in the translating/rotating sphere problem, as well as the squirmer example, that when $h/\epsilon \approx 10$ or greater, the errors appear to saturate. We interpret this as the dominance of the regularization error over the errors created by the imposed linearization of forces and the spatial discretization. However, when $\epsilon$ is chosen much smaller than $h$, we observe the expected second order convergence in $h$, i.e. the regularization error is much smaller than the error due to the force linearization/discretization.

On the other hand, there is a lower bound to how small one can make $\epsilon$ relative to $h$. As we show in \ref{regchoice} of the Appendix, some analytic integrals cannot be evaluated due to the limits of finite precision arithmetic. This leads to a practical lower bound of $\epsilon^2 > \text{eps}\left(\text{max}_{\triangle_i}(L_i) \right)$ where $\text{max}_{\triangle_i}(L_i)$ is the largest side length of a triangle and $\text{eps}(x)$ gives the distance between a positive floating-point number $x$ and the next largest floating-point number.

This leads us to a subtle point about the analytic expressions presented in this work: while it is true that the single layer integral with the singular Stokes kernel is only weakly singular (i.e. the integrand diverges but \textit{is integrable} over a surface), we cannot simply let $\epsilon \to 0$ in \eqref{eq:stokessurfsolparam} and recover a nondivergent expression. This occurs due to the way we decomposed the integral into subintegrals $T_{m,n,q}$ in \eqref{eq:Tmnq}. It is possible that a subintegral $T_{m,n,q}$ diverges in the limit $\epsilon \to 0$ even though when taken altogether as a single integrand, we have an integrable expression.

We note that in the boundary element methods community, analytic approaches to evaluating the nonregularized single layer integral have made recent progress. Ravnik \cite{ravnik2023} and Bohm \textit{et al.} \cite{bohm2024} independently derived expressions for the single layer integrals in the case that the field point lies on the triangular boundary element using a constant force density. The work presented here differs with respect to our use of a piecewise linear force density and a regularized kernel. Although the regularization introduces some error, we have shown in the examples that this error can be made small enough to be negiligible compared to the discretization error. As a final point, from a practitioner's perspective it may be preferable to have a single expression for the integrals in the case when the field point lies on a triangular element and in the case when it does not, which is not yet possible in the nonregularized framework. Using a regularized framework, one does not have to make this distinction.

Future work on this topic could investigate the benefits of higher degree polynomial functions for the forces (e.g. quadratic, cubic), as well as the use of curved triangles to better approximate smooth surfaces. Another potential direction is in developing the method for regularized Stokeslet kernels whose associated blob functions have higher-order moment conditions, as these can lead to greater accuracy without significant increase in computational cost.

\newpage
\bibliographystyle{unsrt}
\bibliography{references}
\appendix
\section{Appendix} \label{appendix3}

\subsection{Derivation of $A_{m,n,q}$ and $B_{m,n,q}$} \label{ABderiv}
Refer to Figure \ref{fig:triangle} for accompanying picture.
\begin{equation}
    \begin{split}
      A_{m,n,q}&= \int_{\partial \Omega_p} \alpha^m(\theta) \beta^n(\theta) \ R^{-q} \ \whvec{n}(\theta) \cdot \whvec e_1 \ d\theta \\
               &= \int_{S_1} \alpha^m(\theta) \beta^n(\theta) \ R^{-q} \ \whvec{n}(\theta) \cdot \whvec e_1 \ d\theta + \int_{S_2} \alpha^m(\theta) \beta^n(\theta) \ R^{-q} \ \whvec{n}(\theta) \cdot \whvec e_1 \ d\theta \\
      &+ \int_{S_3} \alpha^m(\theta) \beta^n(\theta) \ R^{-q} \ \whvec{n}(\theta) \cdot \whvec e_1 \ d\theta\\
      &= \int_{S_1} \alpha^m(\theta) \beta^n(\theta) \ R^{-q} \ \left ((-\whvec e_2)\cdot \whvec e_1 \right) \ d\theta + \int_{S_2} \alpha^m(\theta) \beta^n(\theta) \ R^{-q} \ \left (\whvec{e_1} \cdot \whvec e_1 \right ) \ d\theta \\
      &+ \int_{S_3} \alpha^m(\theta) \beta^n(\theta) \ R^{-q} \left ( \whvec d \cdot \whvec e_1 \right) \ d\theta\\  
      &= \int_0^1 \theta^n R^{-q}(\vec x_f, \vec y_{S_2}(\theta)) \ d\theta - \frac{\sqrt{2}}{\sqrt{2}}  \int_0^1 (1-\theta)^{m+n} R^{-q}(\vec x_f, \vec y_{S_3}(\theta)) \ d \theta \\
      &= \int_0^1 \theta^n R^{-q}(\vec x_f, \vec y_{S_2}(\theta)) d\theta - \int_0^1 (1-\theta)^{m+n} R^{-q}(\vec x_f, \vec y_{S_3}(\theta)) \ d\theta \\
      &= S^{\whvec{e_2}}_{n,q} - \sum_{k=0}^{m+n}\binom{m+n}{k}(-1)^k S^{\whvec{d}}_{k,q}
    \end{split}
    \label{eq:Amnq}
  \end{equation}where $S^{\whvec e_2}_{n,q}$ and $S^{\whvec d}_{k,q}$ can be evaluated using \eqref{eq:S0m1}, \eqref{eq:S0p1}, and \eqref{eq:Srec}. The calculation for $B_{m,n,q}$ is similar but there are two cases, depending if $n=0$ or $n>0$:

  \begin{equation}
    \begin{split}
      B_{m,n,q}&= \int_{\partial \Omega_p} \alpha^m(\theta) \beta^n(\theta) \ R^{-q} \ \whvec{n}(\theta) \cdot \whvec e_2 \ d\theta \\
               &= \int_{S_1} \alpha^m(\theta) \beta^n(\theta) \ R^{-q} \ \whvec{n}(\theta) \cdot \whvec e_2 \ d\theta + \int_{S_2} \alpha^m(\theta) \beta^n(\theta) \ R^{-q} \ \whvec{n}(\theta) \cdot \whvec e_2 \ d\theta \\
      &+ \int_{S_3} \alpha^m(\theta) \beta^n(\theta) \ R^{-q} \ \whvec{n}(\theta) \cdot \whvec e_2 \ d\theta\\
      &= \int_{S_1} \alpha^m(\theta) \beta^n(\theta) \ R^{-q} \ \left ((-\whvec e_2)\cdot \whvec e_2 \right) \ d\theta + \int_{S_2} \alpha^m(\theta) \beta^n(\theta) \ R^{-q} \ \left (\whvec{e_1} \cdot \whvec e_2 \right ) \ d\theta \\
      &+ \int_{S_3} \alpha^m(\theta) \beta^n(\theta) \ R^{-q} \ \left (\whvec d \cdot \whvec e_2 \right ) \ d\theta \\
      &= -1 \int_0^1 \theta^m \beta^n(\theta) R^{-q}(\vec x_f, \vec y_{S_1}(\theta)) \ d\theta + \frac{\sqrt{2}}{\sqrt{2}} \int_0^1 (1-\theta)^{m+n}R^{-q}(\vec x_f, \vec y_{S_3}(\theta)) \ d\theta \\
      &= -1 \int_0^1 \theta^m \beta^n(\theta) R^{-q}(\vec x_f, \vec y_{S_1}(\theta)) \ d\theta + \int_0^1 (1-\theta)^{m+n}R^{-q}(\vec x_f, \vec y_{S_3}(\theta)) \ d\theta \\
      &= \begin{cases}
        \displaystyle - S^{\whvec{e_1}}_{m,q} + \sum_{k=0}^{m} \binom{m}{k}(-1)^{k}S^{\whvec{d}}_{k,q} & n=0 \\
        \displaystyle \sum_{k=0}^{m+n} \binom{m+n}{k}(-1)^k S^{\whvec{d}}_{k,q} & n>0
      \end{cases}
    \end{split}
    \label{eq:Bmnq}
  \end{equation}The reason for the different formulas when $n$ differs is that $\beta \equiv 0$ along the segment $S_1$. When $n=0$, $\beta^n=1$. When $n>0$, $\beta^n=0$.

\subsection{Evaluation of integral \eqref{eq:S1int}} \label{T003details}
The integral used in the computation for $T_{0,0,3}$ is 
\[  \int_0^1 \frac{1}{\sqrt{(\alpha +P)^2+Q^2} \left(\sqrt{(\alpha +P)^2+Q^2}+\gamma/L\right)} \ d \alpha \]The evaluation of this integral requires several substitutions which are described below.

\subsubsection*{First substitution}

Substitute $Qy=\alpha + P$ so $d\alpha = Q \ dy$. Then we have

\begin{equation}
  \begin{split}
      \label{eq:sub1} &\int_{P/Q}^{(1+P)/Q}
\frac{1}{\sqrt{(Qy)^2+Q^2}(\sqrt{(Qy)^2+Q^2}+\gamma/L)} Q \ dy \\
=\frac{1}{Q}&\int_{P/Q}^{(1+P)/Q}
\frac{1}{\sqrt{y^2+1}(\sqrt{y^2+1}+\gamma/(LQ))} \ dy
  \end{split}
\end{equation}
\subsubsection*{Second substitution}

Substitute $t=\sqrt{y^2+1}$. Here is where two cases must be
considered. \\ If $-1<P<0$, then the dummy variable $y$ in
\eqref{eq:sub1} changes sign from negative to positive over the
interval of integration. When $y$ is negative (on the interval from
$P/Q$ to 0), we have $y=-\sqrt{t^2-1}$ and when $y$ is positive (on
the interval from 0 to $(1+P)/Q$), we have $y=\sqrt{t^2-1}$. The
respective differentials are then $dy= -\frac{t dt}{\sqrt{t^2-1}}$ and
$dy=\frac{t dt}{\sqrt{t^2-1}}$.

In this case, \eqref{eq:sub1} becomes

\begin{equation}
  \begin{split}
      \label{eq:sub2neg} &\frac{1}{Q}\int_{P/Q}^{0}
\frac{1}{\sqrt{y^2+1}(\sqrt{y^2+1}+\gamma/(LQ))} \ dy \\ +&\frac{1}{Q}
\int_{0}^{(1+P)/Q} \frac{1}{\sqrt{y^2+1}(\sqrt{y^2+1}+\gamma/(LQ))} \ dy \\ &=
-\frac{1}{Q}\int_{t_1}^{1}\frac{1}{(t+\gamma/(LQ))\sqrt{t^2-1}} \ dt \\ &+
\frac{1}{Q}\int_{1}^{t_2}\frac{1}{(t+\gamma/(LQ))\sqrt{t^2-1}} \ dt
  \end{split}
\end{equation}On the other hand, if $P \geq 0$ or $P \leq -1$,
\eqref{eq:sub1} becomes
\begin{equation}
  \label{eq:sub2}
  \begin{split} \pm \frac{1}{Q}&\int_{t_1}^{t_2}
\frac{1}{t(t+\gamma/(LQ))}\frac{t}{\sqrt{t^2-1}} \ dt \\ \pm
\frac{1}{Q}&\int_{t_1}^{t_2}\frac{1}{(t+\gamma/(LQ))\sqrt{t^2-1}} \ dt
  \end{split}
\end{equation}where $t_1=\sqrt{(P/Q)^2+1}$,
$t_2=\sqrt{(1+P)^2/Q^2+1}$, and the sign is positive if $P \geq 0$ and
negative if $P \leq -1$.

Now that these cases have been taken care of, the rest of the
substitutions will be the same for all cases. We will follow the case
where $P \geq 0$ but the results for the others are similar. The only
things that change on a case by case basis are the limits of
integration.
\subsubsection*{Third substitution} Let $t=\sec \phi$ so that $dt=\sec
\phi \tan \phi \ d\phi$ and $\phi_1=\sec^{-1}(t_1)$,
$\phi_2=\sec^{-1}(t_2)$, and $\sqrt{t^2-1}=\sqrt{\sec^2 \phi -1} =
\sqrt{\tan^2 \phi}=\tan \phi$. Substitute into \eqref{eq:sub2}:

\begin{equation}
  \label{eq:sub3}
  \begin{split} \frac{1}{Q}&\int_{\phi_1}^{\phi_2} \frac{1}{(\sec \phi
+ \gamma/(LQ))\tan \phi} \sec \phi \tan \phi \ d\phi \\
=\frac{1}{Q}&\int_{\phi_1}^{\phi_2} \frac{\sec \phi}{\sec \phi + \gamma/(LQ)}
\ d \phi \\ =\frac{1}{Q}&\int_{\phi_1}^{\phi_2}\frac{1}{1+\gamma/(LQ) \cos
\phi} \ d\phi
  \end{split}
\end{equation}

\subsubsection*{Fourth substitution} Finally, use the tangent half-angle
substitution, $r=\tan(\phi/2)$. Then $d\phi=\frac{2}{1+r^2}dr$,
$r_1=\tan(\frac{\phi_1}{2})$, $r_2=\tan(\frac {\phi_2}{2})$, and $\cos
\phi = \frac{1-r^2}{1+r^2}$. Substitute into \eqref{eq:sub3} to get
\begin{equation}
  \label{eq:sub4}
  \begin{split}
&\frac{1}{Q}\int_{r_1}^{r_2}\frac{1}{1+\gamma/(LQ)\frac{1-r^2}{1+r^2}}
\frac{2}{1+r^2} \ dr \\
&=\frac{2}{Q}\int_{r_1}^{r_2}\frac{1}{1+r^2+\gamma/(LQ)(1-r^2)} \ dr \\
&=\frac{2}{Q}\int_{r_1}^{r_2}\frac{1}{1+\gamma/(LQ)+(1-\gamma/(LQ))r^2} \ dr \\
&=\frac{2}{Q(1+\gamma/(LQ))}\int_{r_1}^{r_2}\frac{1}{1+\frac{1-\gamma/(LQ)}{1+\gamma/(LQ)}r^2}
\ dr \\ &=\frac{2}{Q(1+\gamma/(LQ))} \sqrt{\frac{1+\gamma/(LQ)}{1-\gamma/(LQ)}} \tan^{-1}\left
( r\sqrt{\frac{1-\gamma/(LQ)}{1+\gamma/(LQ)}} \right ) \Bigg \rvert_{r_1}^{r_2}
  \end{split}
\end{equation}

\subsection{Choice of regularization $\epsilon$} \label{regchoice}
As noted in the conclusion, there is a lower bound to how small one can make $\epsilon$ without running into floating point precision issues.

\begin{figure}[ht]
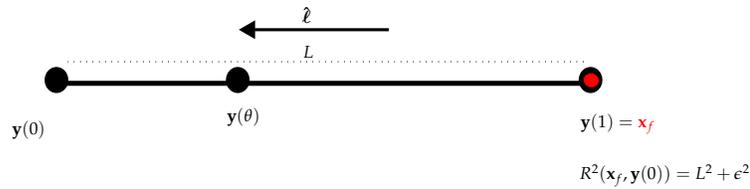

    \centering
    \scalebox{0.6}{
    \def\svgwidth{\columnwidth}
    
}
    \caption[Numerical precision issue for line integrals with small $\epsilon$]{When the regularization, $\epsilon$, is chosen too small relative to $L$, issues can arise in the evaluation of the line integrals \eqref{eq:S0m1} or \eqref{eq:S0p1}.}
    \label{fig:precision}
\end{figure}

For particular values of $\epsilon$, the formulas \eqref{eq:S0m1} and \eqref{eq:S0p1} may lead to singular results. Suppose that the field point, $\vec x_f$, is the same as $\vec y(1)$, the endpoint of the segment shown in Figure \ref{fig:precision}. Then, $\vec x(0) \cdot \bell = (\vec x_f - \vec y(0)) \cdot \bell =  -L$ and $R(\vec x_f, \vec y(0) ) = \sqrt{L^2+\epsilon^2}$. Examining the formula for \eqref{eq:S0m1}, we notice that when evaluating at $\theta = 0$, we have

\[\log \left [ \vec x(0) \cdot \bell + R(\vec x_f,\vec y(0))\right ]\]If the distance from $L$ to the next largest floating point precision number is less than $\epsilon^2$, then $R(\vec x_f, \vec y(0))\approx L+\epsilon^2/(2L)$ will evaluate to $L$ in floating point precision arithmetic, assuming that $L$ is larger than $\epsilon^2$ (this will be the case in practice since $L$ is the length of a triangle side.) The argument of the logarithm is then zero, which is undefined.

A similar issue occurs in \eqref{eq:S0p1}, where the argument of the $\text{arctanh}$ function evaluates to $-1$ in floating precision arithmetic. This is also undefined. These cases give a practical lower bound on how small one can make $\epsilon$. Since the distance between a positive floating point number $x$ and the next largest floating point number increases as $x$ gets larger, the lower bound for $\epsilon$ depends on the length of the largest side of a triangle. For a triangle with maximum side length $L=\text{max}_{\triangle_i}(L_i)$, $\epsilon^2$ must be as large as the distance between $\text{max}_{\triangle_i}(L_i)$ and the next largest floating point number to evaluate the line integrals. 
 
\end{document}

%% file: figures/segment.pdf_tex
\begingroup%
  \makeatletter%
  \providecommand\color[2][]{%
    \errmessage{(Inkscape) Color is used for the text in Inkscape, but the package 'color.sty' is not loaded}%
    \renewcommand\color[2][]{}%
  }%
  \providecommand\transparent[1]{%
    \errmessage{(Inkscape) Transparency is used (non-zero) for the text in Inkscape, but the package 'transparent.sty' is not loaded}%
    \renewcommand\transparent[1]{}%
  }%
  \providecommand\rotatebox[2]{#2}%
  \newcommand*\fsize{\dimexpr\f@size pt\relax}%
  \newcommand*\lineheight[1]{\fontsize{\fsize}{#1\fsize}\selectfont}%
  \ifx\svgwidth\undefined%
    \setlength{\unitlength}{600.07462629bp}%
    \ifx\svgscale\undefined%
      \relax%
    \else%
      \setlength{\unitlength}{\unitlength * \real{\svgscale}}%
    \fi%
  \else%
    \setlength{\unitlength}{\svgwidth}%
  \fi%
  \global\let\svgwidth\undefined%
  \global\let\svgscale\undefined%
  \makeatother%
  \begin{picture}(1,0.27231059)%
    \lineheight{1}%
    \setlength\tabcolsep{0pt}%
    \put(0,0){\includegraphics[width=\unitlength,page=1]{segment.pdf}}%
    \put(0.42536605,0.2402072){\color[rgb]{0,0,0}\makebox(0,0)[lt]{\lineheight{1.25}\smash{\begin{tabular}[t]{l}$\bell$\end{tabular}}}}%
    \put(0.83364861,0.08752877){\color[rgb]{0,0,0}\makebox(0,0)[lt]{\lineheight{1.25}\smash{\begin{tabular}[t]{l}$\vec y(1)$\end{tabular}}}}%
    \put(0.00005659,0.07759279){\color[rgb]{0,0,0}\makebox(0,0)[lt]{\lineheight{1.25}\smash{\begin{tabular}[t]{l}$\vec y(0)$\end{tabular}}}}%
    \put(0,0){\includegraphics[width=\unitlength,page=2]{segment.pdf}}%
    \put(0.42655104,0.18979862){\color[rgb]{0,0,0}\makebox(0,0)[lt]{\lineheight{1.25}\smash{\begin{tabular}[t]{l}$L$\end{tabular}}}}%
    \put(0,0){\includegraphics[width=\unitlength,page=3]{segment.pdf}}%
    \put(0.31503548,0.0962209){\color[rgb]{0,0,0}\makebox(0,0)[lt]{\lineheight{1.25}\smash{\begin{tabular}[t]{l}$\vec y(\theta)$\end{tabular}}}}%
    \put(0,0){\includegraphics[width=\unitlength,page=4]{segment.pdf}}%
  \end{picture}%
\endgroup%

%% file: figures/triPerspective.pdf_tex
\begingroup%
  \makeatletter%
  \providecommand\color[2][]{%
    \errmessage{(Inkscape) Color is used for the text in Inkscape, but the package 'color.sty' is not loaded}%
    \renewcommand\color[2][]{}%
  }%
  \providecommand\transparent[1]{%
    \errmessage{(Inkscape) Transparency is used (non-zero) for the text in Inkscape, but the package 'transparent.sty' is not loaded}%
    \renewcommand\transparent[1]{}%
  }%
  \providecommand\rotatebox[2]{#2}%
  \newcommand*\fsize{\dimexpr\f@size pt\relax}%
  \newcommand*\lineheight[1]{\fontsize{\fsize}{#1\fsize}\selectfont}%
  \ifx\svgwidth\undefined%
    \setlength{\unitlength}{523.04031012bp}%
    \ifx\svgscale\undefined%
      \relax%
    \else%
      \setlength{\unitlength}{\unitlength * \real{\svgscale}}%
    \fi%
  \else%
    \setlength{\unitlength}{\svgwidth}%
  \fi%
  \global\let\svgwidth\undefined%
  \global\let\svgscale\undefined%
  \makeatother%
  \begin{picture}(1,0.81061907)%
    \lineheight{1}%
    \setlength\tabcolsep{0pt}%
    \put(0,0){\includegraphics[width=\unitlength,page=1]{triPerspective.pdf}}%
    \put(0.8681906,0.19941917){\color[rgb]{0,0,1}\makebox(0,0)[lt]{\lineheight{1.25}\smash{\begin{tabular}[t]{l}$\widehat{\vecg{\nu}}$ \end{tabular}}}}%
    \put(0.73691709,0.01742534){\color[rgb]{0,0,1}\makebox(0,0)[lt]{\lineheight{1.25}\smash{\begin{tabular}[t]{l}$\whvec{t}$\end{tabular}}}}%
    \put(0,0){\includegraphics[width=\unitlength,page=2]{triPerspective.pdf}}%
    \put(0.15240921,0.2154323){\color[rgb]{1,0,1}\makebox(0,0)[lt]{\lineheight{1.25}\smash{\begin{tabular}[t]{l}$\vec y_0 - L \alpha \whvec v$\end{tabular}}}}%
    \put(0.27203097,0.29167385){\color[rgb]{1,0,1}\makebox(0,0)[lt]{\lineheight{1.25}\smash{\begin{tabular}[t]{l}$-H \beta \whvec n$\end{tabular}}}}%
    \put(0.80709427,0.78251043){\color[rgb]{0,0,0}\makebox(0,0)[lt]{\lineheight{1.25}\smash{\begin{tabular}[t]{l}$\vec x_f$\end{tabular}}}}%
    \put(0,0){\includegraphics[width=\unitlength,page=3]{triPerspective.pdf}}%
    \put(0.4561561,0.32223076){\color[rgb]{1,0,1}\makebox(0,0)[lt]{\lineheight{1.25}\smash{\begin{tabular}[t]{l}$\vec y(\alpha,\beta)$\end{tabular}}}}%
    \put(-0.00030694,0.22353876){\color[rgb]{0,0,0}\makebox(0,0)[lt]{\lineheight{1.25}\smash{\begin{tabular}[t]{l}$\vec y_0$\end{tabular}}}}%
    \put(0.37811599,0.41997321){\color[rgb]{0,0,0}\makebox(0,0)[lt]{\lineheight{1.25}\smash{\begin{tabular}[t]{l}$\vec y_2$\end{tabular}}}}%
    \put(0.71032136,0.21631695){\color[rgb]{0,0,0}\makebox(0,0)[lt]{\lineheight{1.25}\smash{\begin{tabular}[t]{l}$\vec y_1$\end{tabular}}}}%
    \put(0.40122592,0.22064997){\color[rgb]{0,0,0}\makebox(0,0)[lt]{\lineheight{1.25}\smash{\begin{tabular}[t]{l}$S_1$\end{tabular}}}}%
    \put(0.59621607,0.34053273){\color[rgb]{0,0,0}\makebox(0,0)[lt]{\lineheight{1.25}\smash{\begin{tabular}[t]{l}$S_2$\end{tabular}}}}%
    \put(0.18745784,0.35391065){\color[rgb]{0,0,0}\makebox(0,0)[lt]{\lineheight{1.25}\smash{\begin{tabular}[t]{l}$S_3$\end{tabular}}}}%
    \put(0.80095725,0.04387776){\color[rgb]{0,0,0}\makebox(0,0)[lt]{\lineheight{1.25}\smash{\begin{tabular}[t]{l}$\vec P(\vec x_f)$\end{tabular}}}}%
    \put(0,0){\includegraphics[width=\unitlength,page=4]{triPerspective.pdf}}%
  \end{picture}%
\endgroup%

%% file: figures/radialPolarGraph.pdf_tex
\begingroup%
  \makeatletter%
  \providecommand\color[2][]{%
    \errmessage{(Inkscape) Color is used for the text in Inkscape, but the package 'color.sty' is not loaded}%
    \renewcommand\color[2][]{}%
  }%
  \providecommand\transparent[1]{%
    \errmessage{(Inkscape) Transparency is used (non-zero) for the text in Inkscape, but the package 'transparent.sty' is not loaded}%
    \renewcommand\transparent[1]{}%
  }%
  \providecommand\rotatebox[2]{#2}%
  \newcommand*\fsize{\dimexpr\f@size pt\relax}%
  \newcommand*\lineheight[1]{\fontsize{\fsize}{#1\fsize}\selectfont}%
  \ifx\svgwidth\undefined%
    \setlength{\unitlength}{392.15728519bp}%
    \ifx\svgscale\undefined%
      \relax%
    \else%
      \setlength{\unitlength}{\unitlength * \real{\svgscale}}%
    \fi%
  \else%
    \setlength{\unitlength}{\svgwidth}%
  \fi%
  \global\let\svgwidth\undefined%
  \global\let\svgscale\undefined%
  \makeatother%
  \begin{picture}(1,0.64271993)%
    \lineheight{1}%
    \setlength\tabcolsep{0pt}%
    \put(0,0){\includegraphics[width=\unitlength,page=1]{radialPolarGraph.pdf}}%
    \put(0.30240493,0.40148776){\color[rgb]{0,0,0}\transparent{0.99207598}\makebox(0,0)[lt]{\lineheight{1.25}\smash{\begin{tabular}[t]{l}$a$\end{tabular}}}}%
    \put(0.67149061,0.47523667){\color[rgb]{0,0,0}\transparent{0.99207598}\makebox(0,0)[lt]{\lineheight{1.25}\smash{\begin{tabular}[t]{l}$r=\lvert \vec x \rvert-a$\end{tabular}}}}%
    \put(0.71041369,0.61242544){\color[rgb]{0,0,0}\transparent{0.99207598}\makebox(0,0)[lt]{\lineheight{1.25}\smash{\begin{tabular}[t]{l}$\vec x$\end{tabular}}}}%
    \put(0,0){\includegraphics[width=\unitlength,page=2]{radialPolarGraph.pdf}}%
    \put(0.34709813,0.33107073){\color[rgb]{0,0,0}\transparent{0.99207598}\makebox(0,0)[lt]{\lineheight{1.25}\smash{\begin{tabular}[t]{l}$\phi$\end{tabular}}}}%
    \put(0,0){\includegraphics[width=\unitlength,page=3]{radialPolarGraph.pdf}}%
    \put(0.02067028,0.39195342){\color[rgb]{0,0,0}\transparent{0.99207598}\makebox(0,0)[lt]{\lineheight{1.25}\smash{\begin{tabular}[t]{l}$\vec U$\end{tabular}}}}%
  \end{picture}%
\endgroup%

%% file: figures/poiseuilleRectangle.pdf_tex
\begingroup%
  \makeatletter%
  \providecommand\color[2][]{%
    \errmessage{(Inkscape) Color is used for the text in Inkscape, but the package 'color.sty' is not loaded}%
    \renewcommand\color[2][]{}%
  }%
  \providecommand\transparent[1]{%
    \errmessage{(Inkscape) Transparency is used (non-zero) for the text in Inkscape, but the package 'transparent.sty' is not loaded}%
    \renewcommand\transparent[1]{}%
  }%
  \providecommand\rotatebox[2]{#2}%
  \newcommand*\fsize{\dimexpr\f@size pt\relax}%
  \newcommand*\lineheight[1]{\fontsize{\fsize}{#1\fsize}\selectfont}%
  \ifx\svgwidth\undefined%
    \setlength{\unitlength}{680.34519466bp}%
    \ifx\svgscale\undefined%
      \relax%
    \else%
      \setlength{\unitlength}{\unitlength * \real{\svgscale}}%
    \fi%
  \else%
    \setlength{\unitlength}{\svgwidth}%
  \fi%
  \global\let\svgwidth\undefined%
  \global\let\svgscale\undefined%
  \makeatother%
  \begin{picture}(1,0.37813674)%
    \lineheight{1}%
    \setlength\tabcolsep{0pt}%
    \put(0,0){\includegraphics[width=\unitlength,page=1]{poiseuilleRectangle.pdf}}%
    \put(0.22214446,0.03773231){\color[rgb]{0,0,0}\makebox(0,0)[lt]{\lineheight{1.25}\smash{\begin{tabular}[t]{l}$y=-a$\end{tabular}}}}%
    \put(0.26879821,0.13204903){\color[rgb]{0,0,0}\makebox(0,0)[lt]{\lineheight{1.25}\smash{\begin{tabular}[t]{l}$y=a$\end{tabular}}}}%
    \put(0.56202975,0.08170073){\color[rgb]{0,0,0}\makebox(0,0)[lt]{\lineheight{1.25}\smash{\begin{tabular}[t]{l}$z=-b$\end{tabular}}}}%
    \put(0.55967427,0.26150735){\color[rgb]{0,0,0}\makebox(0,0)[lt]{\lineheight{1.25}\smash{\begin{tabular}[t]{l}$z=b$\end{tabular}}}}%
    \put(0.39182136,0.17303416){\color[rgb]{0,0,0}\makebox(0,0)[lt]{\lineheight{1.25}\smash{\begin{tabular}[t]{l}$2s$\end{tabular}}}}%
    \put(0.44400086,0.12727014){\color[rgb]{0,0,0}\makebox(0,0)[lt]{\lineheight{1.25}\smash{\begin{tabular}[t]{l}$2s$\end{tabular}}}}%
    \put(0,0){\includegraphics[width=\unitlength,page=2]{poiseuilleRectangle.pdf}}%
    \put(0.9760042,0.30764124){\color[rgb]{0,0,0}\makebox(0,0)[lt]{\lineheight{1.25}\smash{\begin{tabular}[t]{l}$x$\end{tabular}}}}%
    \put(0.96029591,0.34972327){\color[rgb]{0,0,0}\makebox(0,0)[lt]{\lineheight{1.25}\smash{\begin{tabular}[t]{l}$y$\end{tabular}}}}%
    \put(0.91368918,0.36567728){\color[rgb]{0,0,0}\makebox(0,0)[lt]{\lineheight{1.25}\smash{\begin{tabular}[t]{l}$z$\end{tabular}}}}%
    \put(-0.04341914,0.32551564){\color[rgb]{0,0,0}\makebox(0,0)[lt]{\lineheight{1.25}\smash{\begin{tabular}[t]{l}$\vec u = (u(y,z),0,0)$\end{tabular}}}}%
    \put(0.08355861,0.04112629){\color[rgb]{0,0,0}\makebox(0,0)[lt]{\lineheight{1.25}\smash{\begin{tabular}[t]{l}$x=-L$\end{tabular}}}}%
    \put(0.80949129,0.04501272){\color[rgb]{0,0,0}\makebox(0,0)[lt]{\lineheight{1.25}\smash{\begin{tabular}[t]{l}$x=L$\end{tabular}}}}%
    \put(0,0){\includegraphics[width=\unitlength,page=3]{poiseuilleRectangle.pdf}}%
  \end{picture}%
\endgroup%

%% file: figures/precisionProblem.pdf_tex
\begingroup%
  \makeatletter%
  \providecommand\color[2][]{%
    \errmessage{(Inkscape) Color is used for the text in Inkscape, but the package 'color.sty' is not loaded}%
    \renewcommand\color[2][]{}%
  }%
  \providecommand\transparent[1]{%
    \errmessage{(Inkscape) Transparency is used (non-zero) for the text in Inkscape, but the package 'transparent.sty' is not loaded}%
    \renewcommand\transparent[1]{}%
  }%
  \providecommand\rotatebox[2]{#2}%
  \newcommand*\fsize{\dimexpr\f@size pt\relax}%
  \newcommand*\lineheight[1]{\fontsize{\fsize}{#1\fsize}\selectfont}%
  \ifx\svgwidth\undefined%
    \setlength{\unitlength}{600.07462629bp}%
    \ifx\svgscale\undefined%
      \relax%
    \else%
      \setlength{\unitlength}{\unitlength * \real{\svgscale}}%
    \fi%
  \else%
    \setlength{\unitlength}{\svgwidth}%
  \fi%
  \global\let\svgwidth\undefined%
  \global\let\svgscale\undefined%
  \makeatother%
  \begin{picture}(1,0.27231059)%
    \lineheight{1}%
    \setlength\tabcolsep{0pt}%
    \put(0,0){\includegraphics[width=\unitlength,page=1]{precisionProblem.pdf}}%
    \put(0.42536605,0.2402072){\color[rgb]{0,0,0}\makebox(0,0)[lt]{\lineheight{1.25}\smash{\begin{tabular}[t]{l}$\bell$\end{tabular}}}}%
    \put(0.83364861,0.08752877){\color[rgb]{0,0,0}\makebox(0,0)[lt]{\lineheight{1.25}\smash{\begin{tabular}[t]{l}$\vec y(1)= \color{red}{\vec {x}_f}$\end{tabular}}}}%
    \put(0.00005659,0.07759279){\color[rgb]{0,0,0}\makebox(0,0)[lt]{\lineheight{1.25}\smash{\begin{tabular}[t]{l}$\vec y(0)$\end{tabular}}}}%
    \put(0.83182182,0.01107084){\color[rgb]{0,0,0}\makebox(0,0)[lt]{\lineheight{1.25}\smash{\begin{tabular}[t]{l}$R^2(\vec x_f, \vec y(0) ) = L^2 + \epsilon^2$\end{tabular}}}}%
    \put(0,0){\includegraphics[width=\unitlength,page=2]{precisionProblem.pdf}}%
    \put(0.42655104,0.18979862){\color[rgb]{0,0,0}\makebox(0,0)[lt]{\lineheight{1.25}\smash{\begin{tabular}[t]{l}$L$\end{tabular}}}}%
    \put(0,0){\includegraphics[width=\unitlength,page=3]{precisionProblem.pdf}}%
    \put(0.31503548,0.0962209){\color[rgb]{0,0,0}\makebox(0,0)[lt]{\lineheight{1.25}\smash{\begin{tabular}[t]{l}$\vec y(\theta)$\end{tabular}}}}%
    \put(0,0){\includegraphics[width=\unitlength,page=4]{precisionProblem.pdf}}%
  \end{picture}%
\endgroup%